\documentclass[11pt]{amsart}
\usepackage{amsmath,amsfonts,amssymb,amsthm,epsfig,cancel}
\usepackage[mathscr]{euscript}
\usepackage{hyperref} 
\usepackage[all]{xy}
\usepackage{color}
\usepackage[normalem]{ulem}
 
\newtheorem{thm}{Theorem}[section]
\newtheorem*{*thm}{Theorem}
\newtheorem{lemma}[thm]{Lemma}
\newtheorem{prop}[thm]{Proposition}
\newtheorem{propdef}[thm]{Proposition-Definition}
\newtheorem{corr}[thm]{Corollary}
\theoremstyle{definition}
\newtheorem{dfn}[thm]{Definition}
\newtheorem{exmple}[thm]{Example}
\newtheorem{exmples}[thm]{Examples}

\theoremstyle{remark}
 
\newtheorem*{rmq}{\textit{Remark}}
\newtheorem{rmk}[thm]{\textit{Remark}}
\renewcommand{\proof}{\noindent\textit{Proof}\/: \,\,}

\newcommand{\C}{{\mathbb{C}}}
\newcommand{\Q}{{\mathbb{Q}}}
\newcommand{\R}{{\mathbb{R}}}

\newcommand{\Z}{{\mathbb{Z}}}

\newcommand{\bH}{{\mathbb{H}}}

\newcommand\germ[1]{{\mathfrak{#1}}}

\newcommand\DD{{\mathcal D}}

\newcommand{\comp}{\raise1pt\hbox{{$\scriptscriptstyle\circ$}}}
\def\lset{\{}  
\def\rset{\}}  
\def\set#1{\lset#1\rset} 
\def\st{\mid}   
\def\sett#1#2{\lset #1 \st #2 \rset}  
\def\del{\partial}

\newcommand\Tr{{}^{\mathsf{T}}\kern-0.9pt} 
\def\id{\mathop{\rm id}\nolimits}
\newcounter{lijstc}
{\end{list}}
{\end{list}}

\newenvironment{diagram}{
\begin{matrix}}{\end{matrix}}

\def\arrow(#1,#2)\dir(#3,#4)\long#5{\put(#1,#2){\vector(#3,#4){#5}}}

\def\mapright#1{\mathop{\vbox{\ialign{
                                ##\crcr
    ${\scriptstyle\hfil\;\;#1\;\;\hfil}$\crcr
 \noalign{\kern2pt\nointerlineskip}
    \rightarrowfill\crcr}}\;}}

\def\mapleft#1{\mathop{\vbox{\ialign{
                                ##\crcr
    ${\scriptstyle\hfil\;\;#1\;\;\hfil}$\crcr
 \noalign{\kern2pt\nointerlineskip}
    \leftarrowfill\crcr}}\;}}

\newcommand\rarrow[3]{\smash{\mathop{\hbox to#3{\rightarrowfill}}\limits
^{\scriptstyle#1}_{\scriptstyle#2}}}
 
\newcommand\larrow[3]{\smash{\mathop{\hbox to#3{\leftarrowfill}}\limits
^{\scriptstyle#1}_{\scriptstyle#2}}}

\newcommand\Vline[3]{\llap{$\scriptstyle #1$}
\left\Vert\vbox to#3{}\right.\rlap{$\scriptstyle #2$}}
\renewcommand\setminus{-}
\def\into{\hookrightarrow}
\def\onto{\twoheadrightarrow}
\def\ii{{\rm i}}

\newcommand\gl[1]{\operatorname{GL}({#1})}
\newcommand\so[1]{\operatorname{SO}({#1})}

\newcommand\ogr[1]{\operatorname{O}({#1})}
\newcommand\ugr[1]{\operatorname{U}({#1})}

\newcommand\smpl[1]{\operatorname{Sp}({#1})}

\newcommand\im{\operatorname{Im}}

\newcommand\Hom{\mathop{\rm Hom}\nolimits}

\def\bF{\mathbb{F}}
\newtheorem{rmks}[thm]{\textit{Remarks}}
\def\clif#1#2{{\mathsf C}^{#1}(#2)}
\def\cl#1{\clif{}{#1}}
\def\gl#1{\text{Gl}(#1)}

\def\spin#1{\text{\rm Spin}(#1)}
\def\sps#1#2{\mathsf{S}^{#1}(#2)}
\def\dirac#1{\cancel{#1}}
\def\ind#1{\text{\rm ind}(#1)}
\def\td#1{\text{\rm td}(#1)}
 
  \def\aut{\text{\rm Aut}}
\def\ch#1{\text{\rm ch}(#1)}
\def\half{ \frac 12}
\def\rot#1{%
\textcolor{red}{%
#1}%
}
\def\half{ \frac 12}

\def\textoplus{{\textstyle\bigoplus}}
\begin{document}

\title[Abelian Varieties associated Riemannian manifolds]{Abelian varieties   and theta functions associated to  compact Riemannian manifolds; constructions inspired by superstring theory.\\
}
\author[S. M\"uller-Stach, C.  Peters and V. Srinivas]{S. M\"uller-Stach, C.  Peters and V. Srinivas\\
Math.  Inst., Johannes Gutenberg Universit\"at Mainz,\\
Institut Fourier, Universit\'e Grenoble I\\
St.-Martin d'H\`eres, France and\\
TIFR, Mumbai, India}
\date{Feb-06, 2012}
\maketitle

\section{Introduction}

\subsection{Motivation from Physics} 
In some forms of superstring theory particular theta functions  come up as  partition functions. The associated abelian varieties come either from certain cohomology groups of the underlying "universe"   or,   in more recent theories (e.g. \cite{wit} \cite{world}), are linked to  their K-groups. 

Expressed  in mathematical terms,     one canonically associates   to the cohomology or the K-theory of an  even dimensional compact spin  manifold    a principally polarized  abelian variety. Moreover, if the dimension is $2 \bmod 8$  a particular  line bundle is singled out whose  first Chern class is  the principal polarization.  This bundle thus has a non-zero section represented by a   theta function which, after suitable normalization, is  indeed the partition function of the underlying theory. 

 There are several  types of superstring theories, e.g. type I which is self-dual and     types IIA and  IIB which are related via $T$-duality. The theories start from a space-time $Y$ which in a first approximation can be taken to be $Y=X\times T$ where $T$ is the time-``axis''\footnote{ $T$ could   be a circle in physical theories but it could even  be a point (``absence of branes'').}  and $X$ is some compact Riemannian manifold. In Type IIA theory the \emph{Ramond-Ramond field}  is a closed differential form $G=G_0+G_2+\cdots$ on $X$ with components of all even degrees while in type IIB $G$ is an odd degree closed differential form on $X$.  Moreover, these forms are \emph{integral} (that is they have integral periods over integral homology cycles). The reason is that  they are Poincar\'e dual to certain submanifolds of $X$ which are the ``world''-part of a brane in $Y$.  Such a field  should be thought of as some configuration  in the theory. The partition function assembles all possible configurations in some generating function which can in turn be used to derive further physical properties of the model.  In type IIA theory this partition function is  of the form $\Theta(0)/\Delta$ where $\Theta$ is some normalized theta function. While $\Delta$ is canonically associated to the Riemannian manifold $X$, this is no longer the case for $\Theta$. Instead, as suggested by Witten in \cite{wit} and later by Moore and Witten in \cite{world} one should  lift the discussion up to K-theory using  the Chern character. But then, in order to make a canonical choice for $\Theta$ one has to assume that the manifold has a spin structure and has dimension $2 \bmod 8$.  The background from physics is collected in \S~\ref{PhysBack}.  It is not necessary for an  understanding of  the rest of the paper, but it purports to explain how  physicist came to the particular jacobians and the normalized theta functions.
 
For algebraic geometers these constructions  may look    a bit esoteric at first sight, the more since they are  phrased in terms foreign to most of them. For instance, one  might ask: \emph{is the construction related to the Weil jacobian? } This  question   was one of the motivations for the present note. Clearly, an answer entails a careful analysis of the construction proposed in \cite{world}.  

\subsection{Mathematical Contents}

The constructions from \cite{world} use in a critical way the index theorems of Atiyah and Singer, a subject not too well known among algebraic geometers.
On the other hand, people well versed in this topic might not have heard about Griffiths'  period domains. So, to make this paper profitable for readers with widely different backgrounds, chunks of theory from several  branches of mathematics have been summarized in  a way adapted to the needs of this paper.

In \S~\ref{RGBack} we review some basic Riemannian geometry, in \S~\ref{sec:KandIndex} K-theory and related index theorems are summarized while in \S~\ref{sec:HodgeTheory} basic constructions concerning Hodge structures are reviewed. The introductory part finishes with a  short motivating  section (\S~\ref{PhysBack}) on quantum field theory.  

The basic construction implicitly used in \cite{world} is really simple and given in \S\ref{sec:LinAlg}.  It  is apparently well known among physicists    but   we  could not find a reference   for it in this precise  form, although a variant is well known in symplectic geometry, cf. \cite[Prop. 2.48 (ii)]{symplec}. In \S~\ref{ssec:HomoGeneous} this  construction has been phrased in terms of homogeneous spaces.  Noteworthy is a diagram given in Th.~\ref{MainResultAsDiagram} which summarizes this. Curiously, this fits very well    with Teichm\"uller theory which describes moduli of  compact Riemann surfaces of a given genus in terms of conformal equivalence classes of metrics. 

Next, the basic construction is applied in several different situations. First of all, in \S~\ref{ssect:FirstEx}  a polarized torus is associated to  even or odd  cohomology of a given compact Riemannian manifold  whose dimension is $2 \bmod 4$.   These constructions  can be twisted by certain automorphisms of the cohomology. These can be exploited to see that  the construction   generalizes the construction of the Weil intermediate jacobian. So, for polarized complex algebraic manifolds this twisted version is a canonical choice.
 
 The question arises if some other extra structure on a Riemannian manifold  in a similar manner leads to  a canonical choice of abelian variety. The idea is that the extra structure should be such that it comes with   a natural differential operator whose index can be calculated from a unit  in the cohomology ring. The action of this unit  then defines the canonical twist. The main example is a spin-structure where we have the associated Dirac operator. By the  Atiyah-Singer index theorem   its index is calculated using  the $\hat A$-genus which   provides the   unit in the rational cohomology ring.   In \S~\ref{ssect:SecondEx} the reader finds the details. In this case  the abelian variety is principally polarized. This abelian variety  is exactly the one from \cite{world}  
\\
As a parenthesis, it should be noted that in loc. cit. no proof is offered that the polarization on this abelian variety  is  indeed principal.        This is true and we show  this \emph{crucial fact} by  reducing it to   an old result \cite{ah3} on normalized multipliers.    

 As suggested by the Teichm\"uller approach we propose as a moduli space for a given compact smooth manifold the space of conformal metrics on it. The construction then gives  various  period maps associated to the manifold. See  Theorem~\ref{CohAsModuli} and  Theorem~\ref{MainThmForKGroups}.

Next, in \S~\ref{sec:HS} we show how in the context of abstract Hodge structures the construction of the Weil jacobian fits in this framework. This applies to  \emph{odd} weight. However, there is an apparently new construction related to even weight which also leads to a polarized abelian variety. By means of an  example with holomorphically varying Hodge structures, we show that  the new abelian variety   varies  in general   non-trivially with parameters. However, as the example shows, as in the case of  Weil jacobians, the dependence is  in general  neither    holomorphic nor anti-holomorphic. In \S~\ref{ssec:CohKaehler} this abstract construction is applied to the cohomology of K\"ahler manifolds (Theorem~\ref{CohKaeherAndModuli}). The reader should contrast this with the tori obtained  using the cohomology of  general Riemannian manifolds.  See  Theorem~\ref{MainThmForKGroups}.

The note ends  with \S~\ref{MainResult} where, after    a short digression on normalized theta functions a mathematical formulation is offered of the pertinent  results of \cite[\S 3]{world}.  
Isolating the line bundle from the numerical equivalence class  of the principal polarization uses in a crucial way   some  constructions  from  real K-theory.   These are quite subtle and have been placed in appendix B.  Noteworthy in this appendix   is  a version  of the Thom-isomorphism in real K-theory (=Theorem~\ref{ThomIsoReal}) which extends the one found in the literature  for spin-manifolds whose dimension is divisible by $8$. This generalization can be extracted from \cite{at2} but is not explicitly stated there.  We state and prove it in the appendix  since this form of the Thom-isomorphism theorem is used in a crucial way in  \cite{wit} to find the ``right'' $\theta$-function.

\subsection*{Acknowledgement}  \begin{small}  Our thanks go to Stefan Weinzierl for his helpful explanation of some of the physical aspects of the theory and SFB/TRR 45 for financial support.
 \end{small}
 
\subsection*{Notation} For any topological space $X$ we let $H^k(X)$ be the $\Q$-vector space of the $k$-th singular cohomology of $X$ with \emph{rational} coefficients. The sublattice $H^k(X)_\Z$ of   images  of  integral cohomology classes is canonically isomorphic to integral cohomology modulo torsion. Furthermore we put
 \[
 H^k(X)_\R:= H^k(X)\otimes_\Q \R,\qquad H^k(X)_\C:= H^k(X)\otimes_\Q \C
 \]
If $X$ is a manifold, we let $T(X)$ be its tangent bundle and $T(X)^\vee$ its dual, the co-tangent bundle of $X$. If $X$ is a smooth (i.e., $C^\infty$-)manifold we denote the vector space of smooth differential $p$-forms on $X$ by $A^p(X)$.
 \\
 We put
 \[
 H^+= \textoplus_{k\in\Z} H^{2k}(X),\quad H^-=\textoplus_{k\in\Z} H^{2k+1}(X)
\]
so that we have a canonical $\Z/2\Z$-graded $\Q$-algebra
\[
H^*(X)= H^+(X)\oplus H^-(X).
\]
For an oriented compact connected manifold $X$ the intersection form 
\[
 H^*(X)\times H^*(X) \to \Q,\quad \langle x,y \rangle= \int_X x\wedge y
\]
is a perfect pairing (Poincar\'e duality) and if moreover $\dim X$ is even, it induces perfect pairings 
on $H^\pm(X)$. On $H^-(X)$ it induces a skew form
 {
\[
\omega^-: H^-(X)\times H^-(X)\to \Q,\quad \omega^-(x,y)=\int_X x\wedge y,
\] 
but on $H^+(X)$ a symmetric form. We replace it by a skew form as follows.}
\\
The splitting   $H^+(X)= H^{4*}(X)\oplus H^{4*+2}(X)$ defines the  involution $\iota$ with the first subspace as $+1$-eigenspace and the second as $(-1)$-eigenspace. Suppose next that $\dim X\equiv 2 \bmod 4$. Then the form  
\[
\omega^{ {+}} : H^+(X)\times H^+(X) \to \Q,\quad \omega^{ {+}}(x,y)=\int_X x \wedge \iota(y) 
\]
is a perfect skew pairing.
Both pairings $\omega^\pm$ restrict to unimodular integral skew pairings  on $H^{\pm}(X)_\Z$.
We set
\[
\omega:=\omega^++\omega^-.
\]
 
 \section{Background from Riemannian Geometry}\label{RGBack}

\subsection{The Hodge Metric} Let $(X,g)$ be a compact $d$-dimensional Riemannian manifold. The metric $g$ induces a metric on the co-tangent bundle and its $k$-th exterior powers, the bundles of $k$-forms. This is a fibre-wise metric.

Next, assume that $X$ has an orientation. Then the Hodge $*$ operator can be defined as follows. Let $\set{e_1,\dots,e_d}$ be an oriented orthonormal basis for $T(X)^\vee$. Let $I=\set{i_1,\dots,i_k}$ be an ordered subset of the ordered set $[n]:=\set{1,\dots,n}$. Then the $e_I = e_{i_1}\wedge \cdots\wedge e_{i_k} \in \Lambda^kT(X)^\vee$ give an orthormal basis for $\Lambda^kT(X)^\vee$ and one defines 
\begin{eqnarray*}
*: \Lambda^kT(X)^\vee & \to & \Lambda^{d-k}T(X)^\vee \\
e_I & \mapsto &e_{[n]\setminus I}.
\end{eqnarray*}
This defines a linear   operator $*: A^k(X) \mapsto A^{d-k}(X)$. We need the following property (see \cite[4.10(6)]{warn}):
\begin{lemma} Let $\dim X=d$ be even. Then $*^2= (-1)^k$ on $H^k(X)_\R$. In particular, on odd cohomology it defines a complex structure.
\end{lemma}

Note that for all $I$ the element $e_I\wedge * e_I=   *1$ is the volume form on $X$ and the fibrewise metric on $\Lambda^kT(X)^\vee $ defines a metric on $A^k(X)$ given by
\[
\langle \alpha ,\beta \rangle=\int_X \alpha \wedge *\beta.
\]
By definition, the operator $d^*:A^k(X) \to A^{k-1}(X)$ is the   adjoint of $d$ with respect to these metrics, i.e.,
\[
\langle \alpha ,d\beta\rangle = \langle d^*\alpha,\beta\rangle   \text{ \rm   for all  forms }\alpha,\beta.
\]
A form $G$ is \emph{co-closed}, respectively \emph{co-exact} if $d^*G=0$, respectively  $G=d^* H$ for some $(d+1)$-form $H$.  A form is \emph{harmonic}  if it is closed and co-closed. 

The \emph{Hodge decomposition theorem} \cite[Chapter 6]{warn} states that there is an orthogonal decomposition
\[
\begin{matrix}
A^p(X)&=& \mathsf{Har}^p(X) & \oplus & d\mathsf{A}^{p-1}(X) & \oplus& d^*\mathsf{A}^{p+1}(X)\\
\Vline{}{}{3ex} && \Vline{}{}{3ex} &&   \Vline{}{}{3ex}&& \Vline{}{}{3ex} \\
[\text{$p$-forms}]   & &  [\text{harmonic $p$-forms}] && [\text{exact $p$-forms}]& &[\text{co-exact $p$-forms}] .
\end{matrix}
\]
Moreover, the space of harmonic forms is \emph{finite dimensional} and every De Rham cohomology class has a unique representing harmonic form.
From harmonic theory it  also follows  (see loc. cit.) that $d\oplus d^*$ induces  a  self-adjoint  operator on
$\mathsf{Har}(X) ^\perp\subset \mathsf{A}(X)$.  

Finally we remark that the metric on $A^k(X)$ when restricted to the subspace of harmonic forms defines a metric on $H^k(X)_\R$, the \emph{Hodge metric}. 
\begin{dfn} \label{HodgeMetric} Denote the cohomology class of a closed form $\alpha$ by $[\alpha]$. 
The \emph{Hodge metric} on $H^* (X)_\R$ is the metric associated to $g$   defined by
\[
b^{(g)}([a],[b]):= \int_X \alpha  \wedge *\beta, \qquad a=[\alpha], \, b=[\beta]
\]
and  $\alpha,\beta$ are the unique harmonic forms in the classes $a,b $ respectively.  \end{dfn}

\subsection{K\"ahler Manifolds}\label{ssec:KahlerMets}
Assume $X$ is a compact complex manifold of dimension $d$.  Then $X$ is a real $(2d)$-dimensional manifold with an almost complex structure $J$.  
Decompose the hermitian metric $h$ on $T(X)$ into real and imaginary parts
\[
h(x,y)= g(x,y)+ \ii \omega(x,y).
\]
Then  $h$ is called \emph{K\"ahler}, if the real  non-degenerate skew-symmetric form $\omega$  is closed.  Since $g(x,y)=\omega(x,Jy)$ and hence is determined by $\omega$ we sometimes call the pair $(X,\omega)$ a K\"ahler manifold. The Riemannian metric is then written $g_\omega$.

A manifold admitting a K\"ahler metric is called a \emph{K\"ahler manifold}.
Examples include projective space with the Fubini-Study metric and projective manifolds with metric induced from the Fubini-Study metric.

\subsection{Spin Manifolds}
For this section the reader may consult \cite[Chapters 3,4]{heat}.
For spin-groups consult Appendix A.  {Recall that  for any vector bundle $E$ on a compact topological space $X$ we have functorially behaving Stiefel-Whitney classes $w_k(E)\in H^k(X; \Z/2\Z )$ in cohomology with $\Z/2\Z$-coeffcients. If $E$ happens to be a complex bundle we also have the Chern classes $c_k\in H^{2k}(X;\Z)$ in integral cohomology. They are related: $w_{2k}(E)$ is the reduction modulo two of $c_k(E)$, i.e., the image of $c_k(E)$ under the coefficient morphism $H^*(X;\Z)\to H^*(X;\Z/2\Z)$.}
\begin{dfn} (1) \label{dfn:Spin} Let $X$ be an oriented Riemannian manifold. An even rank vector bundle $E$ has a  \emph{spin$^c$-structure} if for some \emph{integral class} $w\in H^2(X;\Z)$ one has $w\equiv w_2(E) \bmod 2$ 
\\
(2) If $T(X)$ has a spin$^c$-structure, we say that   $X$ \emph{has spin$^c$-structure}. Note that in particular this implies that $\dim X$ should be even.
\\
(3)   A \emph{spin structure} on a $d$-dimensional manifold $X$ is a $\spin d$-principal bundle $\spin X$ on $X$ such that
\[
T(X)^\vee \simeq \spin X\times_{\spin d} \R^d.
\]
Any manifold having a spin structure is said to be a \emph{spin manifold}.
\end{dfn}
\begin{rmks} \label{spinexample}  (1) In particular, $E$ has a spin$^c$-structure if $w_2(E)=0$. It is well known that this  latter  condition for $E=T(X)$,  i.e., $w_2(X)=0$ is equivalent to $X$  being  spin. In particular, any spin manifold of even dimension has a spin$^c$-structure.  It is well known that the number of in-equivalent spin-structures equals the rank of $H^1(X,\bF_2)\simeq H_1(X,\bF_2)$ (by the universal coefficient theorem). In particular, there are at most finitely many such structures and if $X$ is simply connected there exists at most  one spin structure.   \\
(2) The preceding  remark makes it easy to find examples of spin-manifolds:  compact Riemann surfaces, (real or complex) tori, complex K3-surfaces, and,  more generally  any complex manifold whose canonical bundle is a square.
\end{rmks}
 In Appendix A we recall the  notion of a Clifford algebra.  This notion can be globalized to  the framework of vector bundles on a compact Riemannian manifold $(X,g)$. The cotangent bundle $T(X)^\vee $ is a metric bundle. The Clifford algebras $\cl {T_x^\vee}$, $x\in X$ glue together to give the \emph{Clifford algebra} $\cl X$.  The Riemannian structure defines a unique metric connection \footnote{with respect to the metric $g$; it is a $g$-connection}  on $T(X)$ (and on $T(X)^\vee $) without torsion, the \emph{Levi-Civita connection}. Both $g$ and this connection extend to   the entire exterior algebra $\Lambda^*  T(X)^\vee $ and, using the isomorphism  of  Lemma~\ref{WedgeIsCliff} produces a $g$-connection and an induced Levi-Civita  connection $ \nabla^{\rm LC}$ on the Clifford-algebra. The Clifford algebra has a self-adjoint Clifford action. More generally  one  defines:
\begin{dfn} A \emph{Clifford bundle} is a triple $( W  ,h,\nabla)$ consisting of a  $\Z_2$-graded complex  $\cl X $-module $ W = W ^+\oplus  W^-$ equipped with a hermitian metric $h$ and an  $h$-connection $\nabla$ such that
\\
i) The Clifford action  on the module $ W $, denoted $c:  \cl X \to \aut( W )$, is graded, $ W ^+$ and $ W ^-$ are mutually $h$-orthogonal  and the action is self-adjoint with respect to $h$, i.e.,  
\[
h(c(\alpha) \sigma,\tau )+ h(\sigma,c(\alpha) \tau)=0,
\]
for all differentiable   sections $\sigma,\tau $ of $ W $ and differential  $1$-forms  $\alpha$\, \footnote{Recall that $\alpha(x)\in T_x^\vee\subset \cl{T_x^\vee}$.}
;
\\
ii) The connection is compatible with the Levi-Civita connection in the sense that for any  local vector field  $\xi$ one has
\[
\nabla_\xi (c(\alpha) \sigma) = c(\nabla^{\rm LC}_\xi  \alpha )  \sigma+ c(\alpha) (\nabla_\xi\sigma),
\]
for all differentiable   sections  $\sigma$ of $ W $,  and differentiable $1$-forms $\alpha$.
\end{dfn}
The \emph{Dirac operator} associated to a Clifford bundle $ W $ is  a first order differential operator on the space of sections of $ W $  which is  defined as follows:
\[
  \Gamma(W) \mapright{\nabla} \Gamma(W\otimes_\C  (T(X)^\vee  _\C))    \mapright{c} \Gamma(W) .
\]
It sends sections in $W^\pm$ to sections in $W^\mp$.
\begin{exmples}   {1. The original Dirac operator is Diracs answer as how to find  a square root of the positive Laplacian $\nabla:=-\sum_{j=1}^4\del^2/\del^2 x_j$ on classical space-time $\R^4$ equipped with the Lorentz metric. The Clifford algebra is known to be the algebra of $(2\times 2)$-matrices with coefficients in the quaternions $\bH=\R^4=\C\oplus\C$ with real basis $\set{1,\mathbf{i},\mathbf{j},\mathbf{k}}$ and multiplication rules $\mathbf{i}^2=\mathbf{j}^2=\mathbf{k}^2=-1$, $\mathbf{i}\mathbf{j}=-\mathbf{j}\mathbf{i}=\mathbf{k}$. A complex basis is $\set{1,\mathbf{j}}$.
The Clifford bundle $V$ in question is $\bH^2$. The Dirac operator then is defined to be
\[
\cancel D:= \sum_{k=1}^4 \gamma_k {\del \over \del x_k},
\]
where the $\gamma_k$ are certain $(4\times 4)$-matrices with complex coefficients involving the Pauli-matrices
\[
\sigma_1=\begin{pmatrix} 0 &1  \\ 1 & 0 \end{pmatrix},
\,\sigma_2=\begin{pmatrix} 0&-\mathbf{i}\\\mathbf{i} &0 \end{pmatrix},
\,\sigma_3=\begin{pmatrix} 1 & 0 \\ 0 & -1 \end{pmatrix}.
\]
Precisely, one has
\[
\gamma_1= \mathbf{i}\cdot \begin{pmatrix} \mathbf{1}_2 & 0 \\ 0 & -\mathbf{1}_2 \end{pmatrix}
\, , \gamma_k=\begin{pmatrix} 0& \sigma_k  \\ -\sigma_k&0\end{pmatrix},\; k=2,3,4
\]
for the Euclidean case. In the Lorentz case one removes the factor $\mathbf{i}$ in $\gamma_1$.}
2. The bundle $ \cl X _\C= \cl X \otimes\C$ with its hermitian extension of $g$ and Levi-Civita connection $\nabla^{\rm LC}$  is a Clifford bundle with Dirac operator $d+d^*$. \\
3. Assume $X$ is spin manifold of even dimension $d$ and let   $\mathsf{S} $ be the irreducible complex   $\spin d$-spinor  space (see \eqref{eqn:SSpace} in Appendix A) and form  the \emph{spinor bundle}
\[
\mathbf{\mathcal{S}}= \spin X \times_{\spin d} \mathsf{S}. 
\]
It is a Clifford bundle  when equipped with  the Levi-Civita connection $\nabla^{\rm LC}$ coming from restricting the usual Levi-Civita connection to the subbundle $\mathbf{\mathcal{S}}$ of the complexified Clifford algebra $\cl X_\C$.
The associated Dirac operator is called \emph{the Dirac operator}   of the spin-manifold $X$.
\\
Let $ (E,h) $ be any hermitian vector bundle on $X$ with   { an $h$-metric} connection $\nabla$. The twisted bundle  $W= E\otimes \mathbf{\mathcal{S}}$  has  a product hermitian structure and a natural product connection which is compatible with this metric. All Clifford bundles are of this form. The associated Dirac-operator  $\cancel{D}_E$ is called the \emph{Dirac operator with coefficients in \  $E$}.
\end{exmples}

\section{Summary of $K$-theory and Index Theory}\label{sec:KandIndex}
The reader may consult the excellent introduction \cite{at0}.   For a solid introduction to  index theorems  consult the appendices  \cite[\S 24-26]{hir} to  Hirzebruch's classic.
\subsection{K-theory for Complex Bundles}
We let $X$ be a topological space  and we let $K(X)$ be the Grothendieck group of \emph{complex} vector bundles on $X$.   This is by definition the free $\Z$-module generated by the isomorphism classes of complex vector bundles modulo the relations $E\oplus F -E-F$.  It can be seen to be generated by \emph{virtual bundles}, i.e., differences of the form $E-F$, where $E$ and $F$ are any two vector bundles.
The tensor product on vector bundles is compatible with these relations and so $K(X)$ becomes a ring. If $f:X\to Y$ is continuous, pull back of bundles induces  a ring homomorphism $f^*:K(Y)\to K(X)$.

The \emph{suspension} $SX$ is obtained from the product $S^1\times X$ by identifying all points in the subspace $\set{\mathbf{1}}\times X$ where $\mathbf{1}=(1,0)\in S^1\subset  \R^2$;  $n$-fold iterated suspension is denoted $S^nX$. One defines
\[
K^{-n}X:= K(S^nX).
\]
Bott's periodicity theorem \cite{bo} can be stated as $K^{-2}(X)\simeq K(X)$ which makes it possible to define $K^nX$ for all integers $n$. There are natural pairings $K^{n}(X)\times K^{m}(X)\to K^{n+m}(X)$ compatible with the Bott periodicity  making $\oplus_{n\in \Z}K ^n(X)$ into a graded ring.
In view of Bott's theorem the essential part of this ring is
\[
K^*(X)=K(X)\oplus K^1(X)
\]
with $\Z_2$-grading. The cohomology-ring can also be given a $\Z_2$-grading as $H^*(X)=H^+\oplus H^-(X)$.
The Chern character  (see e.g. formula \eqref{eqn:ChChar} below) gives  $\Z_2$-graded isomorphisms  \cite{ah2}
\begin{equation}\label{eqn:ChernIsIso}
\ 
\text{\rm ch} : K^*(X)\otimes\Q  \mapright{\sim}  H^*(X).
\end{equation}
One can   also define relative $K$  groups $K^n(X,Y)$ where  $Y$ is a subset of $X$ and these fit in exact sequences as for     ordinary cohomology. One
 important fact is the $K$-theoretic version of the \emph{Thom isomorphism theorem}:
 Let $B(E)$, respectively  $S(E)$ be the unit disk-bundle, unit sphere bundle associated to a hermitian  vector bundle $(E,g)$ of rank $r$. Then
 \[
 K^*(B(E),S(E)) \simeq K^{*}(X).
 \]

 \subsection{The Index Theorem}
 
 Let $X$ be a differentiable manifold, $E,F$ two hermitian  vector bundles, $\mathsf{D}:\Gamma(E)\to \Gamma(F)$ a differential operator with adjoint $\mathsf{D}^*$.
 Recall that the \emph{index} is given by
 \begin{equation}\label{eqn:AnINdex}
 \ind {\mathsf D}:= \dim \ker \mathsf D -\dim \ker \mathsf D^*.
 \end{equation}
Let us write for brevity
\[
BX:=B(T(X)^\vee,\qquad SX:=S(TX)^\vee.
\]
The symbol of ${\mathsf D}$ defines an element $\sigma(\mathsf{D})\in K(BX,SX)$ whose Chern character lands in $H^*(BX,SX)\simeq H^*_c(T(X)^\vee )$. Let $\pi: T(X)^\vee \to X$ be the natural projection. Since $T(X)^\vee $ is a symplectic  manifold and $g$ a Riemannian metric on $X$, the bundle  $T(T(X)^\vee )$ has a natural complex structure (Prop.~\ref{AlgProp}) such that $\pi^*TX\otimes\C\simeq T(T(X)^\vee )$ and there is a Todd class 
\[
\td{T(X)^\vee |BX}\in H^*(BX).
\]
Since $H^*(BX,SX)$ is a $H^*(BX)$-module  the following formula makes sense; it defines the  topological index
\begin{equation}\label{eqn:TopInd}
\text{\rm ind}_\tau(\mathsf{D}):= \int_{T(X)^\vee } \ch {\sigma(\mathsf {D})}\cdot  \td {T(X)^\vee |BX}.
\end{equation}
One has:
\begin{thm}[\protect{\cite{AS1}}] Let $X$ be a compact differentiable manifold and ${\mathsf D}$ an elliptic differential operator between complex vector bundles on $X$. Then the topological index \eqref{eqn:TopInd} equals the (analytical) index \eqref{eqn:AnINdex}. In particular, it is an integer.
\end{thm}
The $K$-theoretic extension comes from the remark that  the right hand side of \eqref{eqn:TopInd} makes sense  if we replace  $\sigma(\mathsf {D})$ by any element $d\in K(BX,SX)$.  The topological index for such an element is then defined by 
\[
 \text{\rm ind}_\tau(d)= \int_{T(X)^\vee } \ch  d  \cdot \td {T(X)^\vee |BX}.
 \]
It can be shown that it also makes sense to speak of an \emph{analytic index} $\ind d$ for such elements and that it is an integer:
\begin{thm}[\protect{\cite{AS2}}]  \label{asindexthm}
For a compact differentiable manifold $X$, the two homomorphisms
\[
\text{\rm ind}, \text{\rm ind}_\tau : K(BX,SX) \to \Q
\]
coincide and hence take values in $\Z$.
\end{thm}
This version has a relative form for  differentiable locally trivial fibrations  $f:X\to T$.  The starting observation is the fact that $K(\mbox{\rm point})=\Z$ so that the integer $\ind {\mathsf D}$  is just  the $K$-theoretic  difference of the vector spaces $\ker \mathsf D-\ker \mathsf D^*$. For a family over $T$ this pointwise construction  gives a difference of complex bundles on $T$ and hence an element of  $K(T)$.  For the topological index one  has to replace $T(X)^\vee $ by the relative cotangent bundle $T(X/T)^\vee $ and one gets:
\begin{thm}[\protect{\cite{AS3}}]
For a differentiable  family  $f:X\to T$ of compact differentiable manifolds, the two homomorphisms
\[
\text{\rm ind}, \text{\rm ind}_\tau : K(B(X/T),S(X/T)) \to  K(T)\otimes\Q
\]
coincide and hence take values in $K(T)$.
\end{thm}

\subsection{The Index Theorem for the Dirac Operator}
We start with a few preliminaries on genera. See  \cite[Ch. 1]{hir}. Start
with  any formal power series  $p(z)= 1+ p_1z+p_2z^2+\cdots \in 1+\Z[z]$ whose $m$-th order truncation has  a  formal factorization  
\[
1+p_1 z+\cdots+p_m z^m  = (1+\beta_1 z)\cdots (1+\beta_m  {z}).
\]
Next we  explain a certain formal procedure which uses the $\beta_j$ and a  second formal powerseries of the form
\[
q(z)=1+q_1z+q_2z^2+\cdots \in \Q  [[ z]] 
\]
as an input and whose output is the so called \emph{$q$-series for $p$}.
For some fixed  $ m$   write down the $m$-fold product series
\begin{eqnarray*}
q(\beta_1z)\cdot q(\beta_2z)\cdot\, \dots\,\cdot  q(\beta_m z)& = &
1+Q_1(p_1)z+Q_2(p_1,p_2)z^2+\cdots\\
&=&1+q_1p_1 z+\left[ (-2q_2+q_1^2)p_2- q_1p_1^2\right]z^2+\cdots ,
\end{eqnarray*}
where by definition  the $Q_j$ are the coefficients of $z^j$. These  turn out to be   universal polynomials of total degree $j$ in the first $j$ ``variables'' $p_1,\dots,p_j$ with coefficients expressible in the coefficients
of the formal series $q(z)$. To find these, one calculates  successively, setting  $m=1$, $m=2$ etc.  
The  \emph{$q$-series for $p$} is the resulting formal series
\[
1+Q_1(p_1)z+Q_2(p_1,p_2)z^2+\cdots\in  \Q [ p_1,p_2,\dots][[z]].
\]
The corresponding $q$-genus is obtained by setting $z=1$.
The particular choice 
\[
 q(z)= \frac{ \half z}{ \sinh \half z}= 1- \frac{1}{2^2}\cdot \frac 1 6 z + \frac{1}{2^4}\cdot\frac{7}{360} z^2+\cdots
\]
defines  the $\hat A$-series
\[
\hat A(z,p_1,p_2,\dots))= 1-\frac{1}{2^2}\cdot\frac1 6 p_1z +\frac{1}{2^4}\left(-\frac {1}{90} p_2+ \frac{7}{360} p_1^2\right) z^2+\cdots\, .
\] 
Next we recall that for any   vector bundle $F$ the \emph{Pontryagin classes} are obtained from   the Chern classes as follows: 
\[
p_i(F):= (-1)^i c_{2i}(F\otimes\C)\in H^{4i}(X), i=1,\dots, m=\text{\rm rank} (F).
\]
One associates to these the \emph{Pontryagin polynomial} $p(F):= 1+p_1(F) z+\cdots+p_m(F) z^m $.  The Pontryagin  classes of $X$ are those of $T(X)$.   If we now substitute     $p_i(F)$ for $p_i$ in the  $\hat A$-series,  truncate the series at order $m=\text{\rm rank} (F)$ and  set $z=1$ we obtain the \emph{$\hat A$-genus} for $F$: 
\[
\hat A(F):= \hat A(1,p_1(F),\dots,p_m(F))\in H^{4*}(X;\Q).
\]
In particular 
\[
\hat A(X):= \hat A(T(X)).
\]
For any complex vector bundle $F$ one also has  the Chern character $\ch F$, defined as follows. Write formally $1+c_1(F)x+\cdots c_m(F) x^m=(1+\gamma_1x)\cdots(1+\gamma_mx)$ and evaluate
\begin{equation}\label{eqn:ChChar}
\ch F=\sum e^{\gamma_i}= m+ c_1(F)  + \half (c_1^2(F)-c_2(F))+\cdots \in H^{2*}(X;\Q).
\end{equation}
 Now we have
\begin{thm}[Atiyah-Singer index theorem, \cite{AS1}] \label{ATSIT} Let $X$ be a manifold with a spin structure, $E$ a  complex  bundle   on $X$.  Let $\cancel D_E$ be the Dirac operator  with coefficients in $E$.\footnote{It also makes sense to speak of the Dirac operator with values in a virtual bundle $[E]\in K(X)$ and the same formula for its index holds. It is this version which we need below.} We have
\[
\ind {\cancel D_E}= \int_X \hat A(X) \ch E.
\]
\end{thm}
This is a special case of the general index theorem, Theorem~\ref{asindexthm}.

\begin{rmk} \label{IndexForSpinc}   There is another important  index theorem for  $X$    a compact complex   manifold  carrying a complex vector bundle  $E$. To explain it, let $\Lambda^{p,q}(E)$ be the \emph{bundle} of complex $E$-valued forms  of type $(p,q)$ and denote its sections (the corresponding \emph{forms}) by $A^{p,q}(X)$. Decompose   the $d$-operator  in the usual way as  
\[
d=\del+\bar \del: A^{p,q}(E) \to A^{p+1,q}(E) \oplus A^{p,q+1}(E).
\]
Choose a hermitian metric $h$ on $X$  and $h'$ on $E$  and let $\del^*$, $\bar\del^*$ be the $(h,h')$-adjoints  of $\del$ and $\bar\del$. The bundle 
\[
\textoplus_q \Lambda^{0,q}(E) =\textoplus_k \Lambda^{0,2k+1} (E)  \oplus    \textoplus_k \Lambda^{0,2k}(E)
\]
is a Clifford bundle with $\bar\del+\bar\del^*$ acting on its global sections as  Dirac operator. Applying the index theorem yields   the Hirzebruch-Riemann-Roch theorem 
\[
\int_X\text{\rm td}(X) \ch E = \ind{ \bar\del+\bar\del^*}
\]
and so the left hand side is an integer as well.
\end{rmk} 
Let us see what the above index theorems  gives if $E$ is real. Since $c_i(\bar E)= (-1)^ic_i(E)$ we see from \eqref{eqn:ChChar} that the complex conjugation fixes the terms in  $\ch E$ of degree divisible by $4$ while it acts as minus the identity on the other terms. So the Atiyah-Singer index theorem implies
\begin{corr} \label{Skew} \footnote{It is exactly at this point that  passing   to real operators and real K-theory yields  finer invariants. See Theorem ~\ref{RealIndexThm}.} If $\dim X\equiv 2 \mod 4$  one has $\ind {\cancel D_{\bar E}}=-\ind {\cancel D_E}$.  In particular, if $E$ is real, the index of the Dirac    operator with values in $E$. vanishes.\end{corr}

\section{Background From Hodge Theory}\label{sec:HodgeTheory}
We recall a number of general facts about polarized Hodge structures  \cite[Chapter 2.1]{mhb}.
 
\subsection{Hodge Structures} Recall that a \emph{rational Hodge structure of weight $k$} 
consists of a rational vector space $W$ together with a  decomposition
\[
W_\C:= W\otimes\C=\oplus_{p+q=k} W^{p,q}
\]
such that $W^{q,p}=\overline{W^{p,q}}$. Its \emph{Hodge numbers} are 
\[
h^{p,q}=\dim_\C W^{p,q}
\]
Sometimes we write just $W$ for the structure. It can alternatively   be  described as an algebraic  representation   of the  (real) matrix   group 
\[
\mathbf{G}(\R)= \sett{ s(a,b):=  \begin{pmatrix} a& -b\\
b& a
\end{pmatrix} } 
{a,b\in \R,\, a^2+b^2\not=0} ,
\]
say $h: \mathbf{G}(\R) \to \aut{W_\R}$, $W_\R=W\otimes\R$ with the extra two properties of being  defined over the rationals  and for which   $h(t)= t^k$,  encoding the weight. The Hodge decomposition can then be found back as follows. First identify  $s(a,b)$ with the complex number $z=a+\ii b$  establishing an isomorphism $\mathbf{G}(\R)\simeq  \C^*$ of real algebraic groups.  Then $W^{p,q}$ is the subspace of $W_\C=W\otimes\C$ on which $z\in \C^*$  acts as $z^p\bar z^q$. 

The simplest example making use of this description is the one-dimensional \emph{Tate Hodge structure}   $\Q(-k)$, the one-dimensional representation of $\mathbf{G}(\R)$ given by the character $h(z)= |z|^{k}$. It is pure of type $(k,k)$.

The  Weil operator $C_W= h(\ii)$ is the real operator acting as $\ii^{p-q}$ on $W^{p,q}$. Note that $C_W^2=(-1)^k$ and so defines a complex structure if the weight is odd.
In general the Weil operator alone does not suffice to determine the Hodge structure since it does not determine the $W^{p,q}$. But we have
\begin{lemma} \label{frc} Giving a weight $1$ Hodge structure on $W$ is equivalent to giving a  complex structure on $W_\R:=W\otimes\R$.
\end{lemma}
\proof The space  {$W^{1,0}$, respectively $W^{0,1}$  is the eigenspace of $C_W$ with eigenvalue $\ii$, respectively $(-\ii)$.}  This shows that one must have  $h(a+b\ii)= a\mathbf{1} +b C_W$. Conversely, given a complex structure $J$ on a rational vector space $W$ defines a weight $1$ Hodge structure by defining $h(a+b\ii)= a \mathbf{1} +  bJ$.
\qed\endproof
 
\subsection{Polarizations} \label{ssec:Pol}
A \emph{polarization} for  a weight $k$ Hodge structure $W$ is a non-degenerate $\Q$-valued bilinear form $Q$ on $W$ which is symmetric if $k$ is even and skew-symmetric otherwise and whose complex-linear extension $Q_\C$  satisfies the two Riemann conditions \footnote{Weil demands instead that $Q(x,C_Wx)>0$ for $x\not=0$. The difference changes  the  sign for $Q$  by $(-1)^{k}$.}
\begin{enumerate}
\item $Q_\C(x, y)= 0$ if $x\in W^{p,q}$, $y\in W^{r,s}$ and $(p,q)\not=(s,r)$;
\item $Q_\C(C_W x,\bar x)>0$ if $x\not=0$.
\end{enumerate}

Note that the first condition implies that the Weil operator  {preserves the polarization}: 
$Q(C_Wx,C_Wy)=Q(x,y)$.  Conversely, in weight $1$ we have:
\begin{lemma} \label{fbl}  For a  weight $1$ Hodge structure $W$, the first Riemann bilinear relation  
relative to  a  skew form $Q$ on $W$ is equivalent to  $C_W$ being  $Q$-symplectic.
\end{lemma}
\proof  If $x,y\in W^{1,0}$ one has $Q_\C(x,y)=Q_\C(C_Wx,C_Wy)= -Q_\C(x,y)$.\qed\endproof
\medskip
One can easily see  that the datum of a polarization is equivalent to  giving a   morphism
\[
S: W\otimes W \to \Q(-k),\qquad  S=(2\pi \ii)^{-k}Q 
\]
of weight $(2k)$-Hodge structures for which the bilinear form  defined by $ (x, y) \mapsto (2\pi \ii)^{\rot k} S(C_W x, y)$ is symmetric and positive definite on $W$. 
It follows  that a polarization  induces an isomorphism
\[
\hat Q: W \mapright{\sim} W^\vee (-k)
\]
 of  weight $k$ Hodge structures.
  
If $W=W_\Z\otimes\Q$, i.e.,  $W$ is an \emph{integral Hodge structure}, we speak of an \emph{integral} polarization $Q$ if $Q$ takes integral values on $W_\Z$.   This is inspired by the fact that integral weight one polarized Hodge structures $(W,Q)$ are the same as polarized abelian varieties: to the Hodge structure associate the real torus  $W/W_\Z$  equipped with complex structure induced by the Weil operator $C_W$ and polarization induced by $Q$. The associated polarized abelian variety is then denoted $J(W,Q)$. The polarization is a principal polarization precisely when $Q$ is a unimodular.

The standard example of a polarized Hodge structure comes from the cohomology of a K\"ahler manifold as we now recall briefly.  See e.g. \cite{weil} for details.
\begin{exmple}\label{CohKaehler}
Let $(X,\omega)$ be a compact K\"ahler manifold of dimension $d$.  

Let $L$ be the real operator  on cohomology which is cup product with the class defined by $\omega$.  The weak Lefschetz theorem states that $L^{k} : H^{d-k}(X)_\R  \mapright{\sim} H^{d+k}(X)_\R$. The kernel of $L^{k+1}$ by definition is the \emph{primitive cohomology}
\[
H^{d-k}_{\rm prim}(X)_\R   := \ker \set{L^{k+1}: H^{d-k}(X)_\R  \to H^{d+k+2}(X)_\R }.
\]
By definition $H^k_{\rm prim}(X)_\R =0$ when $k> d$. The Lefschetz decomposition theorem tells us how to build up cohomology from the primitive parts:
\[
H^k(X)_\R  =\bigoplus_{k\ge r} L^m H^{k-2r}_{\rm prim}(X)_\R .
\] 
Write  the primitive decomposition for $x,y\in H^k(X)_\R$ as $x=\sum L^rx_r$, $y=\sum L^sy_s$. Then the
   \emph{Riemann-form} \cite[p.~77]{weil} can be written as:
\begin{equation}\label{eqn:PolarAll}
\begin{array}{lcl}
Q_\omega(\sum_r L^r x_r,\sum_s L^s y_s) &:= & \epsilon_{k}   \sum_r (-1)^r  \mu_r \int_X
 L^{d-k+2r}(x_r\wedge y_r),\\
\hspace{9em}\epsilon_k &:=&(-1)^{\half k (k+1)} , \quad     \mu_r:= \frac{r!}{(d-k+r)!} .
\end{array}
\end{equation}
Riemann's bilinear relations tell us that $Q_\omega$ is symmetric and non-degenerate if $k$ is even and symplectic if $k$ is  {odd}. 

From  Weil's formula \cite[Chap 1, Th. 2]{weil}  for the $*$ operator in terms of the Weil operator $C$ for the Hodge structure on cohomology:
\begin{equation}\label{eqn:Weil}
*(L^r x_r) =   \epsilon_{k}   (-1)^r  \mu_r  L^{d-k+r} Cx_r,\qquad x_r\in H^{k-2r}(X)_{\rm prim} 
\end{equation}
we see that the Hodge metric is given by
\begin{equation}\label{eqn:HodgeComesFromPol}
b^{(g_\omega)}(x,y)=  Q_\omega( x, Cy) ,\quad x,y\in H^k(X).
\end{equation}
If $X$ is a projective manifold, the class $\omega$ can be taken to be rational so that the Lefschetz decomposition as well as the form $Q_\omega$ is rationally defined. 
Remark however that even if we choose for $\omega$ an integral class,   formula \eqref{eqn:PolarAll} shows that the polarization becomes a priori only rational on   $H^k(X)_\Z$.  Not only a  {denomin}ator is introduced but also, even if $x$ is integral, the primitive constituents are in general only rational.    To remedy this, using a basis, one can always find an integer $N$ such that $NQ_\omega$ becomes integral on $H^k(X)_\Z$. Taking $N$ minimal when varying over all possible bases gives an intrinsically defined integer, say $N_\omega$ so that $N_\omega Q_\omega$ becomes integral. 
\end{exmple}

\subsection{Griffiths Domains} We also make use of the Griffiths period domains \cite{periodarticle} which parametrize polarized  real Hodge structures of given weight and Hodge numbers.  

To explain this roughly, let $(W,Q)$ be a given real Hodge structure  of weight $k$. It is characterized by a Hodge flag $F=\set{ F^k\subset F^{k-1}\subset\cdots \subset F^\ell}
$  where $F^k=H^{k,0}$, $F^{k-1}=H^{k,0}+H^{k-1,1}$, etc. If $k=2\ell-1$ (odd case) or $k=2\ell$ (even case) the flag stops at stage $\ell$. The first Riemann condition is an algebraic condition while the second is open in the complex topology. The group $G=\so{W,Q}$ acts transitively on flags satisfying these conditions and the stabilizer of the given flag $F$ is a compact subgroup $H\subset G$ so that the period domain can be written as $D=G/H$.   Below, in \S~\ref{sec:HS}
we recall this in some more detail.

The standard example of a polarized Hodge structure of weight $k$ is the primitive $k$-th cohomology of any K\"ahler variety $X$ as recalled in Example~\ref{CohKaehler}. 
 If moreover $X$ varies in an algebraic family $\set{X_t}_{t\in T}$, the assigment $t\mapsto k$-th primitive cohomology of $X_t$ induces a holomorphic \emph{period map} $p: T \to D$. This map is in general multi-valued because of monodromy.
For details see  e.g.   \cite[Chap. IV.3]{periodbook}.

 \section{Physics background}\label{PhysBack}
 
 {The constructions are in fact inspired by quantum field theory.  Start with a Riemannian
manifold $Y$. Fields on $Y$ are sections in certain vector bundles on $Y$}. These come in three types: scalar fields, fermionic fields and  gauge fields.  The first type is just a function on $Y$,
the second one is a smooth section in a spinor-bundle on $Y$, and the last are Lie-algebra valued forms on $Y$.
Experimentally observed are gauge fields and fermionic fields.
A scalar field has not been observed yet, but physicist are
desperately trying to observe it (the search for the Higgs boson, for which
the new LHC accelarator has been built).


Given the fields,   the physics is deduced  from  a Lagrangian 
density.  {For a free scalar field $\phi$ the Lagrangian reads} (using Einstein summation) 
\[
L = \half \partial_\mu \phi(x) \partial^\mu \phi(x) - \half  \phi(x) \phi(x).
\]
Its integral over $Y$ is the action:
\[
S = \int_Y  \mathrm{d}^nx \, L,
\]
where   $ \mathrm{d}^nx$ is some suitable  $n$-dimensional measure on $Y$, $n= \dim Y $.
The action obviously depends on the field.

 {The manifolds $Y$ occurring in physics are called \text{space-time worlds} and they are to have  a \emph{time component}; moreover,  one assumes that}
there is a Lorentzian metric on $Y$ which is a positive definite metric on world-sheets.
To change it to a positive definite  metric on \emph{all} of $Y$, there is a trick, called \emph{Wick rotation}: one formally replaces the time parameter $t$ by $\ii t$. 
Note that this replaces $S$ by $\ii S$ as well.  {Below we explain why this is a useful trick.}

Once this has been done,  by definition, the partition  function  is   the integral over all fields, where each field   is weighted by $\exp(\ii S)$. For a free scalar field one gets
$Z = \int \mathcal{D}\phi \cdot \exp(\ii S)$.
 {Such an  integral, a path integral, is   mathematically ill-defined.
In physics literature it has  been interpreted in analogy with similar integrals in probability theory and  gets  replaced by an ordinary, but  complicated  integral  over the  so-called configuration space of $Y$. It turns out that this integral can be approximated in a very particular way  up to any  given order in a way described by 
the associated Feynman graphs and Feynman rules for the fields. To actually calculate the integrals in the resulting expansion, one has to  rewrite them  so as to involve   certain contour integrals over paths which are  in the complex plane (thanks to the Wick rotation explained before). Residue theory then makes it possible to calculate these. For several quantum field theories the thus calculated  path integral yields results which  are surprisingly close  to  the experimentally observed values.}

The theory that is important here is inspired by  gauge fields on $4$-manifolds $Y$. These  have
  a nice geometric reformulation for the
action. A gauge field  can be viewed as a
connection  one-form on  a given principal   fibre bundle on $Y$. This is  a Lie algebra valued one-form on 
$Y$. Its  covariant derivative $F=DA$, i.e., the curvature of the fibre bundle is an ordinary $2$-form on $Y$.
The action for the gauge fields can be written in this geometric
language as $S = \int F \wedge *F$,
where $*F$ is the Hodge-star of $F$, another two-form on $Y$.

In string theory one replaces the ordinary $4$-dimensional time-world  $Y$ by some
other variety  of  dimension $10$. There are several types of superstring theories: type I which is self-dual,
and types IIA and IIB which are each others dual. 
In this note we mainly consider  type IIA theories.
For such a theory a  gauge field $F$, an integral closed form,  is replaced by an arbitrary even degree  closed form $G = G_0 + G_2 + ... $ with integral periods,
a \emph{Ramond-Ramond field}
and, in analogy to ordinary gauge theory, the action is
\[
S= \int_Y G\wedge *G.
\]
Next,  one wants to define a partition function of the form 
\begin{equation}\label{eqn:PartFun}
Z = \int \mathcal{D}C\cdot  \exp (\ii S),
\end{equation}
where $\int \mathcal{D}C$  is the path integral over all $C$ with $dC=G$ and where $G$ runs over the even degree closed forms with integral periods. 
In order to make this more precise, one assumes that $Y= X\times T$ where $T$ is the ``time-axis''-
which may or may not be compact and represents time  (usually $T=\R$, but $T=$ a point is also a possibility) and $X=(X,g)$ is a Riemannian manifold of dimension $d$ which is assumed to be compact or at least one on which the Hodge decomposition theorem holds (\S~\ref{RGBack}).  Recall   that $d\oplus d^*$ induces  a  self-adjoint  operator on $\mathsf{Har}(X) ^\perp\subset \mathsf{A}(X)$.
In the theory of Ramond-Ramond fields one is only interested in even degree forms that are already closed. So one writes down the
decomposition $ \mathsf{A}(X)=  \mathsf{A}^+(X )\oplus \mathsf{A}^-(X) $ into even and odd degree forms and to understand the even exact forms, one looks at  the associated Dirac operator: 
   \begin{equation}\label{eqn:formD}
\mathsf{D}:=  \begin{pmatrix} 0 & d \left |_{\mathsf{Har}^-(X) ^\perp}\right. \\
 d^* \left |_{\mathsf{Har}^+(X) ^\perp}\right.    & 0
 \end{pmatrix},
  \end{equation}
   {where $\mathsf{Har}^+(X)$, ($\mathsf{Har}^-(X)$) are the even (odd) degree harmonic forms on $X$.}
 By \cite[\S 9.6]{heat} using $\zeta$-functions there is an exact  way  to define its regularized determinant  $\det\mathsf{D}$ and its square root 
\begin{equation} \label{eqn:RegDet}
\Delta: = \sqrt{\det(\mathsf{D})} 
\end{equation}
is called the   \emph{determinant of the non-zero modes}. In view of the form \eqref{eqn:formD} for the operator    $\mathsf{D}$,  the  determinant of the non-zero modes can be viewed  as the determinant of the operator $d$ on odd degree forms.

Return now to  the integral \eqref{eqn:PartFun}. 
 { Fix $\alpha\in H^{2*}(X)_\Z$ and consider all possible closed forms $G_\alpha$ representing $\alpha$. If $G_\alpha^0$ is the unique harmonic representative in the class $\alpha$ we then can write
\[
G_\alpha= G_\alpha^0+ dC,\quad \langle G_\alpha^0,dC\rangle=0,
\]
where the inner product  has been introduced in~\S 2.1. It follows that the action for the field $G_\alpha$ is a sum
\[
S(G_\alpha)= \underbrace{\langle G_\alpha^0,G_\alpha^0 \rangle}_{S_\alpha^{\rm cl} } + \underbrace{\langle dC,dC\rangle}_{S^{\rm q}},
\]
where the first term depends on $\alpha$, the ``classical'' contribution and the second term, the ``quantum contribution'' does not. All possible $C$ and $\alpha$   together describe all  the field-configurations;  the  partition function becomes
\[
\int \DD C\exp{(\ii S)}=\underbrace{ \left(\int \DD C \exp({\ii S^{\rm q}})\right)}_{Z_{\rm quantum}}\cdot  \underbrace{\sum_{\alpha}  \exp({\ii S_\alpha^{\rm cl}})}_{Z_{\rm classical}}
\]
}%
The relevant 
 calculations  have been carried out in detail in \cite{gauge} \footnote{Actually, in \cite{gauge}  the odd middle degree differential forms  on  a compact manifold of even dimension $d\equiv 2 \bmod 4$ are studied. However, exactly the same calculations  apply in the setting of the article \cite{world}.} and they  can than be summarized as follows:   
\begin{itemize}
\item the classical contribution is of the form
\[ Z_{\rm classical}= \text{Anomalous pre-factor}\cdot \Theta(0) 
\]
 where $\Theta(z)$ is  some classical theta function.
\item  in the total partition function the anomalous pre-factor  {together with  $Z_{\rm quantum}$}  gives   the  factor
$\Delta^{-1}$ where $\Delta$ is the  determinant of non-zero modes \eqref{eqn:RegDet}:
\[
Z= \frac{1}{\Delta} \Theta(0).
\]
\end{itemize}
\begin{rmq}
In loc. cit.   there is given no clue  as to which $\Theta$-function should be used. This is one of the issues  which  the two articles \cite{wit,world} address and which we  want to discuss below in \S~\ref{MainResult}.
\end{rmq}

\section{The basic construction} \label{sec:basic}
\subsection{A Linear Algebra Construction} \label{sec:LinAlg}
To motivate the construction, recall  how one defines  a hermitian metric on a complex  vector space  in terms of real geometry. 
\begin{dfn} Let $V$ be  a real vector space $V$  equipped with a complex structure $J$. A \emph{hermitian metric} on $(V,J)$  is given by a real bilinear form $h:V\times V\to \C$  such that, writing $h$ in real and imaginary parts as $h=b+\ii\omega$ one has
\begin{enumerate}
\item $h$ satisfies $h(Jx,Jy)=h(x,y)$ for all $x,y\in V$;
\item $b$ is a metric and $\omega$ is a  symplectic form, i.e., a non-degenerate skew-symmetric real form;
\item $\omega(x,y)=b(Jx,y)$ for all $x,y\in V$ -- or, equivalently $b(x,y)=\omega(x,Jy)$;
\end{enumerate}
The form $\omega$ is called the metric form  and $b$ the underlying real metric. Condition 3) states that the metric form is uniquely determined by the underlying real metric and the complex structure.\\
 A symplectic  form $\omega$ on $V$ which satisfies the following two weaker conditions
\begin{enumerate}
\item  $\omega(Jx,Jy)=\omega(x,y)$ for all $x,y\in V$;
\item $\omega(x,Jx)>0$ if $x\not=0$
\end{enumerate}
is said to be \emph{tamed by} $J$. 
\end{dfn}
One can ask whether given a (real) metric $b$ and a non-degenerate  skew-symmetric real form $\omega$ 
determine a complex structure so that the two come from a hermitian metric. Note that one should then have 
that   (3) holds. This means that in any case the form $b_{\omega,J}(-,-)=\omega(-,J-)$ should define a 
metric. Below we show that, replacing (3) by this weaker condition, there does exists a unique  complex 
structure such  that $b_{\omega,J}+\ii\omega$ is   hermitian   with respect to $J$. However $b_{\omega,J}$ 
rarely coincides with the original metric $b$. If this happens,  we speak of a coherent pair   $(b,\omega)$ 
(see Definition~\ref{DefCoherent}).{We have the following result.}

 \begin{prop}\label{AlgProp}
Let $(V,b,\omega)$ be a finite dimensional $\R$-vector space equipped with a (positive definite) metric $b$ 
and  a non-degenerate  $\R$-valued skew- symmetric form   $\omega$. There exists a unique complex structure 
$J$ on $V$ such that
\begin{enumerate}
\item $b(Jx,Jy)=b(x,y)$ for all $x,y\in V$;
\item $\omega(Jx,Jy)=\omega(x,y)$ for all $x,y\in V$;
\item the form  $b_{\omega,J}$ defined by  $b_{\omega,J}(x,y):=\omega(x,Jy)$  is { a (positive definite) metric}. \end{enumerate}
\end{prop}
 \proof
Define $A\in \gl V$ by 
\begin{equation}\label{eqn:DefineA}
\omega(x,y)=b(Ax,y).
\end{equation} 
Then $A$ is $b$-skew adjoint:  $A^*=-A$, where $\ast$  means the $b$-adjoint. Hence   $P=A^*A=AA^*=-A^2$ is   self-adjoint and positive definite with respect to $b$. In particular, $V$ has a $b$-orthonormal basis of $P$-eigenvectors so that the matrix of $P$ becomes  diagonal with positive entries, say  $\lambda_i>0$,  on the diagonal. Replacing these by the positive root $\sqrt{\lambda_i}$ defines the root $Q=P^{\half}$ of $P$. Now write
\begin{equation}  \label{eqn:PolDec}
A= Q J,\quad J:=  Q^{-1}A.
\end{equation}
Since $Q^2= P =-A^2$ and $Q^*=Q$ we have
\[
b(Jx,Jy)=b(Q^{-1}A x, Q^{-1}Ay )=- b(Ax,A^{-1}y)=b(x,y)
\]
so that $J$ is $b$-orthogonal. Since    $Q$  is self-adjoint and positive definite, with respect to $b$ this implies that \eqref{eqn:PolDec} is the unique \emph{$b$-polar decomposition} of $A$. Moreover, since $A^*A=AA^*$, one easily deduces that  $J$ and $Q$ (and also $A$ and $Q$) commute. It follows that
\[
J^2= (AQ^{-1})^2= A^2 Q^{-2}= -\id_V
\]
and hence $J$ is a   complex structure. Next,  $\omega(Jx,Jy)= b(AJx,Jy)=b(JAx,Jy)=b(Ax,y)=\omega(x,y)$ since $J$ and $A$ also commute. 
\\
Finally, $\omega(x,Jy)=\omega(Jx,J^2y)=-\omega(Jx,y)=\omega(y,Jx)$ and hence $b_{\omega,J}$ is symmetric. To show that it is positive definite write $\omega(x,Jy)=b(Ax,Jy)=b(QJx,Jy)=b(JQx,Jy=b(Qx,y)=b(Q^{\half}x, Q^{\half}y)$. Then, using  a $b$-orthogonal basis of $Q$-eigenvectors one sees that $b(Qx,x)>0$ if $x\not=0$.
\\
To show uniqueness, note that $A$ is uniquely defined by $b$ and $\omega$ and hence so is its polar decomposition. We only have to see that   $J$ as characterized by 1)--3) gives a polar decomposition $A=RJ$.
First of all $AJ=JA$ as can easily be seen from 1) and 2). But then $R:= AJ^{-1}=J^{-1}A$ is seen to be self-adjoint. From this and  3) it follows that $R$ is positive definite which finishes the proof of uniqueness.
\qed\endproof
 
\begin{corr} \label{IsRiemannForm} Under the hypotheses of \ref{AlgProp} 
 the form $\omega$ is a Riemann-form for $(V,J)$:
\begin{enumerate}
\item $\omega(Jx,Jy)=\omega(x,y)$ for all $x,y\in V$;
\item the form $   \omega(x,Jy)=b_{\omega,J}(x,y)$ is a symmetric $\R$-valued positive definite bilinear form on $V$  for which $J$ is orthogonal. 
\end{enumerate}
\end{corr}
It follows that  the $\C$-valued form
\begin{equation}\label{eqn:IsAPolar}
h(x,y):= \omega(x,Jy)+\ii \omega(x,y)
\end{equation}
is hermitian (with respect to $J$) and positive definite.  Hence:
 
\begin{corr} \label{CruxCor} Suppose   that $V=\Lambda\otimes \R$ for some lattice $\Lambda\subset V$ of maximal rank and that $\omega$ is  {integer valued}. Then  the torus $V /\Lambda$ is an abelian variety with polarization $h$ given by  equation \eqref{eqn:IsAPolar}. This is a principal polarization if $\omega$ is unimodular. 

Conversely, if $(\Lambda,\omega)$ is a lattice equipped with a non-degenerate integral skew form such that $V=\Lambda\otimes\R$ admits a complex  structure $J$ compatible with $\omega$
and such that $\omega(x,Jx)>0$ for $x\not=0$, the form $\omega$ is a Riemann form for the complex torus $V /\Lambda$.  Its complex structure  comes from  the unique   complex  structure of the lemma  with respect to $\omega$ and $b$, where $b(x,y)= b_{\omega,J}(x,y)=\omega(x,Jy)$.
\end{corr}

\begin{dfn} \label{DEFN} The abelian variety just constructed is denoted $J(\Lambda,b,\omega)$. This can equivalently be phrased in terms of Hodge theory: the triple $(\Lambda,b,\omega)$ defines  a unique polarized weight one Hodge structure whose jacobian is $J(\Lambda,b,\omega)$.
\end{dfn}
To stress that  $b\not= b_{\omega,J}$ in general, we recall the notion of conformal equivalence and introduce a new notion.
\begin{dfn}   \label{DefCoherent}  1)   Two  (indefinite) metrics $b,b'$ on a vector space are \emph{conformally equivalent}, if for some    positive constant    $\lambda$ one has $b(x,y)= \lambda b' (x,y) $ for all $x,y\in V$. \footnote{For the usual notion of conformal equivalence of metrics on a manifold in this definition the constant has to be replaced by a function.}We say that $b$ is conformal to $b'$, or $b'$ is conformal to $b$;\\
2) The pair $(b, \omega)$ is  called a \emph{coherent} pair if the metric $b_{\omega, J}$ is conformal to $b$. If $b'$ is conformal to $b$ and $(b,\omega)$ is a coherent pair, then also $(b',\omega)$ is a coherent pair. \\
3) More generally, if we say that  $(b,\omega)$ and $(b',\omega')$ are conformally equivalent if for some positive constants $\lambda,\mu$ we have $b'= \lambda b$ and $\omega'=\mu\omega$.
\end{dfn}
 
\begin{rmks} 1. Observe  that  while conformally equivalent $(b,\omega)$ and $(b',\omega')$ give the same complex structure, the pair  $(b,-\omega)$ gives $-J$. 
\\
2. One can rephrase  Proposition~\ref{AlgProp}  in terms of symplectic geometry as follows:  \textit{A symplectic structure on a finite dimensional euclidean vector space is tamed by a unique complex structure compatible with the metric}.
 \end{rmks}
 We also want to record how coherent pairs behave under the obvious group actions:
 \begin{lemma} \label{Compatible} The group    $\operatorname{Gl}^+(V)$ operates on metrics and symplectic  forms:
 \[
 b_\gamma(x,y):=b(\gamma x,\gamma y),\quad \omega_\gamma(x,y):= \omega(\gamma x,\gamma y),\quad \gamma\in \operatorname{Gl}^+(V).
 \]
If $(b,\omega)$ is coherent, then so is $(b_\gamma,\omega_\gamma)$. The pair $(b_\gamma,\omega)$ is coherent, precisely if $\gamma\in  \operatorname{Aut}(V,\omega)$.
 \end{lemma}

\subsection{Interpretation in Terms of Homogeneous Spaces} \label{ssec:HomoGeneous}
In what follows one  chooses a basis for $V$ to identify $V$ with $\R^{2n}$ and we assume that the symplectic form in this basis is   the standard symplectic form $\omega_o (x,y) = \Tr y \mathbf{J} x$  where 
\[
\mathbf{J}= \begin{pmatrix}
\mathbf{0}_n &   \mathbf{1}_n \\
- \mathbf{1}_n & \mathbf{0}_n\end{pmatrix}. 
\]
With $b$ the   standard metric on $V$, the pair $(b,\omega_o)$ is a coherent couple with complex structure on $\R^n$ given by the matrix $\mathbf{J}$. 
The symplectic group $\smpl {n}$ is the group of $2n\times 2n$-matrices $T$ with $\Tr T \mathbf{J}T =\mathbf{J}$.

The set  $\mathcal{J} ^+ _n$ of complex structures $J$ on  $\R^{2n}$  which preserve  a given orientation   forms  a homogeneous space under conjugation by elements of $\operatorname{Gl}^+(2n,\R)$. This gives an effective group action since  every complex structure on $\R^{2n}$ can be written as a conjugate of the standard complex structure  structure $\mathbf{J}$.  Furthermore,  the isotropy group of the latter is just $\gl{n,\C}$ under the identification
\begin{equation}\label{eqn:Blocks}
  \begin{pmatrix} A & B \\ -B & A  \end{pmatrix}  \mapsto A+\ii B \in \gl{n,\C}.
\end{equation}
Hence
\[
\operatorname{Gl}^+(2n,\R)/ \gl{n,\C}\mapright{\simeq} \mathcal{J} ^+ _n ,\quad [T]  \mapsto   T^{-1} \mathbf{J} T.
\]
The  following is a reformulation of the definition:
\begin{lemma} A  complex structure $J  \in  \mathcal{J} ^+ _n$  is symplectic  if and only if $\omega_o(Jx,Jy)=\omega_o(x,y)$ for all $x,y\in\R^{2n}$. 
\end{lemma} 
While all    complex structures are conjugate to the standard one under the action of $\operatorname{Gl}^+(2n,\R)$, this is no longer the case under the action of $ \smpl{n,\R}$. Indeed, with $\mathbf{D}_k=\text{\rm diag}(\underbrace{1,\dots,1}_{k},\underbrace{-1,\dots,-1}_{n-k}) $, the matrices 
$
  \begin{pmatrix}
\mathbf{0}_n & - \mathbf{D}_k \\
\mathbf{D}_k & \mathbf{0}_n\end{pmatrix} \in   \smpl{n,\R}
$
are in different orbits under the symplectic group. However, \emph{tamed}  almost  complex  structures are conjugate to the standard one by symplectic matrices.  Indeed, pick any basis which at the same time is symplectic and $g_J$-orthonormal. The last condition precisely  means that in this basis $J$ is given by the matrix $\mathbf J$. We in fact have:
\begin{lemma} \label{FromJtob}The subset  $\mathcal{J}^{\omega_o}_n \subset \mathcal{J}_n$ of complex structures  tamed by\footnote{recall that this means in addition to compatibilty $\omega_o(Jx,Jy)=\omega_o(x,y)$ for all $x,y\in  \R^{2n}$ one demands   $\omega_o(x,Jx)>0$ for $x\not=0$.}  $\omega_o$  form  a homogeneous space
\[
 \smpl{n,\R} / \ugr{n} \mapright{\simeq} \mathcal{J} ^{\omega_o}_n   ,\quad [T]  \mapsto T^{-1}\mathbf{J} T . 
\]
Moreover,  the matrix
\[
\Tr T T =   \Tr J \comp \mathbf{J}   = b_{\omega,J}
\]
is positive definite, i.e.,  defines a metric.  In fact $b_{\omega_o,J}$ is the metric associated to the symplectic form $\omega_o$ and complex structure $J$ in accordance with the definition from Proposition~\ref{AlgProp}.
\end{lemma}
\proof
One has 
\begin{eqnarray}
\Tr J \mathbf{J} J&  = &\mathbf {J \label{eqn:IsSymp}}\\
\Tr J \mathbf{J} &&\!\!\text{is symmetric and positive definite.} \label{eqn:IsPos}
\end{eqnarray}
Note that the first condition states that $J\in \smpl{n,\R}$ and it implies that $\Tr J \mathbf{J} $ is symmetric. 
If  $ T\in \smpl{n}$, then $\Tr T =- \mathbf{J}T^{-1} \mathbf{J} $  and hence
\[
\Tr T T  = - \mathbf{J}T^{-1} \mathbf{J} T =  - \mathbf{J}J = \Tr J \comp \mathbf{J}   = b_{\omega,J}
\]
is  positive definite, which proves the second assertion. \qed\endproof
Now note that the collection $\mathsf{Met}_{2n}$ of all metrics   on $\R^{2n}$  form a homogeneous space under $\operatorname{Gl}^+(2n,\R)$ where the group action is by sending $b$ to $\Tr T b  T$ and 
\[
 \operatorname{Gl}^+(2n,\R)/ \so{2n} \mapright{\simeq} \mathsf{Met}_{2n},\quad [T] \mapsto \Tr  T T.
\]
\begin{corr} The inclusion of the symplectic group in the general linear group induces an equivariant inclusion of homogeneous spaces:
\begin{eqnarray}\label{eqn:CanIncl} 
\mathcal{J} ^{\omega_o} _n =   \smpl{n,\R} / \ugr{n} & \into & \operatorname{Gl}^+(2n,\R)/ \so{2n}=\mathsf{Met}_{2n}\\
J= T^{-1}\mathbf{J}  T &\mapsto & \Tr T T =  \Tr J \comp \mathbf{J}=b_{\omega,J}.
\end{eqnarray}
Any other symplectic form $ \omega$ is of the form  $\omega=\Tr \gamma  \omega_o \gamma $ with   $\gamma\in \operatorname{Gl}^+(2n,\R)$. Then $\mathcal{J} ^{ \omega }_n=\gamma^{-1}  \mathcal{J} ^{\omega_o} _n\gamma$ and  the latter gets embedded in the space of metrics by sending $ \omega$ to  the metric $ b_{\omega,J}$ associated to $  \omega$ and the complex structure $J$.  
\end{corr}
Introduce the map
\[ r: \mathsf{Met}_{2n}\to \mathcal{J} ^{\omega_o} _n,\quad b\mapsto J_b,
\]
where $J_b$ is the   unique complex structure compatible with $b$ and the symplectic form $\omega_o$ as given by Prop.~\ref{AlgProp}. We then have
\begin{lemma}  \begin{enumerate}
\item The map $r$ is a retraction for the inclusion \eqref{eqn:CanIncl}. 
\item  Write $J_b= T^{-1} \mathbf{J}  T$ with $ T$ {symplectic}, the fibre   $r^{-1}J_b$ consists of $  \Tr  T \mathcal{K}_n   T$ with
\begin{equation}\label{eqn:Cone}
\mathcal{K}_n =\sett{ Z\in \C^{n\times n}}{Z^*=Z,\, \text{\rm Re}(Z) >0}.
\end{equation} 
In other words,  $r$ exhibits $\mathsf{Met}_{2n}$ as  a fibre bundle over  $\mathcal{J} ^{\omega_o} _n$ with fibres  translates of the $n^2$-dimensional real cone $\mathcal{K}_n$ under the natural action of the  symplectic  group on the space of metrics $\mathsf{Met}_{2n}$.
\end{enumerate} \end{lemma}
\proof (1) Uniqueness  of $J_b$ implies that if one starts  from   a   complex structure   $J= T^{-1}\mathbf{J} T$ with $T$ symplectic, as above,  applying $r$ to    $b_{\omega_o,J}= \Tr T T $ gives back the complex structure $J$. Hence $r$ is indeed a retraction.
\\
(2)  The proof of Prop.~\ref{AlgProp} gives an explicit expression for $r(b)$, $b\in  \mathsf{Met}_{2n}$  a metric.  Indeed, by  \eqref{eqn:DefineA} we have to apply the $b$-polar decomposition  to the matrix $ T:=b^{-1} \mathbf{J}$ and hence  
\begin{equation}\label{eqn:Construction}
r  (b) =   \left( - (b^{-1} \mathbf{J})( b^{-1} \mathbf{J}) \right)^{-\frac 12}  (b^{-1} \mathbf{J}). 
\end{equation} 
The   fibre $r^{-1}J_b$ consists of positive definite symmetric matrices $B$ for which   $J_B=  J_b$.  Formula \eqref{eqn:Construction} for $B$ implies that  $ \mathbf{J} B  J_B= J_B  \mathbf{J} B $
and so, if $J_B=J_b$ one has
\[
[\mathbf{J} B] J_b= J_b [\mathbf{J} B].
\]
This means that $\mathbf{J} B$ is complex linear for the complex structure $J_b$. Conversely, one sees that if the above equation holds, $r(b)=r(B)=J_b$. 

Now perform the change of basis which transforms $J_b$ into $ \mathbf{J}$, i.e.,  with $T$   the change of basis matrix,  $\mathbf{J} B$ gets transformed into
\[
T \mathbf {J} B T^{-1}= \mathbf{J}  \tilde B,\quad \tilde B:= \Tr T^{-1}  B T^{-1} .
\]
The above  equality follows since  $T$ is symplectic.
Hence $\tilde B$ is also a symmetric matrix. Writing $\mathbf{J}  \tilde B$ out in blocks according to the identification \eqref{eqn:Blocks}, one finds that 
$
 \tilde B \in \mathcal{K}_n
$
and so $J_B=J_b$ if and only if $B \in \Tr T  \mathcal{K}_n T$. 
 \qed\endproof

\subsection{Summary}
Suppose now that $V$ is any real vector space of dimension $2n$ and let $\sf{Met}(V)$ be the space of metrics on it and let $(\Lambda^2 V^\vee )^0\subset \Lambda^2 V^\vee $ denote  the open subset of non-degenerate $2$-forms. Recall (Def.~\ref{DefCoherent}) the notion of coherent  pairs on $V$. Let us denote the set of coherent pairs   $\mathsf{CohPair}(V)$.  Recalling also Lemma~\ref{Compatible}, the above discussion  then leads to 
\begin{thm} \label{MainResultAsDiagram}
We have a  commutative  diagram
\[
\hspace{-1em}
 \begin{xy}
 \xymatrix@C=0.65cm{
 \mathsf{Met}(V)\times (\Lambda^2 V^\vee)^0  \ar@<1pc>[d]^{\rho}  
&      \ar@{_{(}->}[l] {\rm Met}(V) \times \{\omega\} \ar@/_2pc/[dd]_(0.3){  r}  \ar[r]^{\hskip -1em  \cong}  & \operatorname{Gl}^+(2n,\R)/ \hskip -1pt \so{2n}  \\
\mathsf{CohPair}(V)  \ar@<1pc>@{_{(}->}[u] & \ar@{_{(}->}[l]    \im (\beta)   \ar@{_{(}->}[u]  \ar[r]^{\hskip -1em  \cong}     &   \im(j)  \ar@{_{(}->}[u] \\
  &  \mathcal{J}^\omega       \ar[u]^\beta_{\cong}
 \ar[r]^{\hskip -1em  \cong} & \smpl{n,\R} /\hskip -1pt \ugr{n} . \ar[u]_{\cong}^j
} 
\end{xy}
\] 
The map $\beta$ associates to the complex structure $J$ tamed by $\omega$ the metric $b_{\omega,J}$ from Lemma~\ref{FromJtob}
while $r(b,\omega)$ is the complex structure  of  Prop.~\ref{AlgProp} and $\rho=b \comp r$ associates to $(b,\omega)$ the coherent  pair $(b_{\omega,J},\omega)$ where $J=r(b,\omega)$. 

The group $\operatorname{Gl}^+(2n,\R)\times \operatorname{Gl}^+(2n,\R)$ operates naturally on $\mathsf{Met}(V)\times (\Lambda^2 V^\vee)^0$.  The diagonal subgroup permutes middle columns for the various symplectic forms while the subgroup $\smpl{n,\R}\times \id $  preserves the middle column for $\omega=\omega_o$, the standard symplectic form.
The  map $r$ is a retraction for the inclusion map defined by $\beta$. Its fibres are translates of a real  $n^2$-dimensional cone isomorphic to \eqref{eqn:Cone}. 
\end{thm}  
 
\begin{rmks} \label{OneDimCase} 1. Note that for $n=1$ giving an orientation is the same as giving a positive $2$-form up to a positive multiple. Such a form is a symplectic form and so
any orientation preserving complex structure is automatically compatible with this symplectic form. The above then  says that giving such  a  complex structure   is equivalent to giving a metric up to a scalar, i.e.,  a class of a conformal metric. 

This is   related to Teichm\"uller theory as follows.  Let us recall briefly the ingredients. One starts with a  compact oriented (real) surface $X$ of genus $g$. For simplicity, assume $g>1$. Any almost complex structure on $X$ is known to be integrable. One lets $\mathcal{J}^\omega(X)$  be the set of all complex structures on $X$ compatible with the orientation. 
This  space does not have the structure of a manifold but is some contractible subset in an infinite dimensional vector space. By the above remarks this space is the same as the space  $\mathsf{Conf}(X)$ for the conformal equivalence classes of metrics on $X$.  Points of $\mathcal{J}^\omega(X)$ can be seen as  equivalence classes of  pairs $(C,f)$ consisting of a genus $g$ curve together with  a diffeomorphism\footnote{classically, one considers homeomorphisms, but this does not matter.} $f:C\to X$, where $(C,[f])$ and $(C',[f']))$ are equivalent if $f^{-1} \comp (f')$ is a biholomorphic map $C'\to C$. The group $\text{Dif}^+(X)$ of orientation preserving diffeomorphisms of $X$ acts on $\mathcal{J}^\omega(X)$. An  orbit  is an isomorphism class  of genus $g$ curves and so the quotient   $\text{Dif}^+(X)\backslash \mathcal{J}^\omega(X)$  is the moduli space   $M_g$  of genus $g$ curves.  The quotient of $\mathcal{J}^\omega(X)$ by the subgroup of those orientation preserving diffeomorphisms that are isotopic to the identity is the \emph{Teichm\"uller space} $\mathcal {T}_g$. This turns out to be a complex manifold (of dimension $3g-3$), in fact it is biholomorphic to some  ball in $\C^{3g-3}$.
By definition, a point in $\mathcal {T}_g$ consists of an equivalence class of a pair $(C,[f])$ of a genus $g$ curve $C$ together with an isotopy class  $[f]$ of an oriented  diffeomorphism $f:C\to X$,  a so-called  \emph{Teichm\"uller structure} on $C$.  Equivalence for pairs is defined as above.  
The isotopy classes  of orientation preserving diffeomorphisms $f:X\to X$ form a group, the \emph{Teichm\"uller group} or \emph{mapping class group},  denoted $\Gamma_g$. It acts on the  Teichm\"uller structure by composition (but does not change $C$) and  $M_g= \Gamma_g\backslash \mathcal{T}_g$. Alternatively, by   the Dehn-Nielsen theorem, the group $\Gamma_g$ is an index 2 subgroup inside the quotient group $ \aut (\pi_1(X))/\set{\text{inner automorphisms}}$. It follows that giving a Teichm\"uller structure $[f]$ on  $C$  is the same as giving  the class of the induced group isomorphism
$f_*: \pi_1(C) \mapright{\sim} \pi_1(X) $ up to inner automorphisms of the target.    If we use   the standard presentation of $\pi_1(X)$ we get the standard symplectic intersection form on the resulting basis of the $1$--cycles so that  the   natural homomorphism $\pi_1(X)\to H_1(X)$ induces a homomorphism $\Gamma_g\to \smpl{g}_\Z$ whose kernel $T_g$ is called the \emph{Torelli group}. 
A \emph{marking} for $C$ is an isometry $  H^1(C)\mapright{\sim} H^1(X).$
The symplectic group  acts on markings by composition.  Over $\mathcal{J}^\omega(X)$ one has a tautological  family of genus $g$ curves. This family descends to Teichm\"uller space and since $\mathcal{T}_g$ is contractible, it  is 
differentiably trivial. Hence the local system of $1$-cohomology groups can be trivialized on Teichm\"uller space and by construction, also on its quotient by the Torelli group. Markings thus globalize over these spaces. Of course this holds also over $\mathcal{J}^\omega(X)$.
 
 {Next note that $J(C)$, the jacobian of $C$ is just obtained  via the construction of   Corr.~\ref{CruxCor} using the   metric  on $H^1(X,\R)$ induced by the conformal class of the metric on $C$ and  the   cup-product pairing on $H^1(X,\Z)$.  Using the global marking, the assigment $C\mapsto J(C)$ defines the period map $\mathsf{Conf}(X) \to \bH_g$ which descends to a holomorphic map 
$p: \mathcal{T}_g \to \bH_g$ and further down to $T_g\backslash \mathcal{T}_g$. Torelli's theorem states that it descends to an injective morphism $M_g \to \smpl{g}_\Z \backslash \bH_g$. 
Summarizing, we have 
\[
\begin{xy}\xymatrix{
\mathsf{Conf}(X)     \ar[r]^{\cong} &  \mathcal{J}^\omega(X)\ar[d] \ar[r] & {\mathbb H}_g\ar@{=}[d]\  \\
& \mathcal{T}_g \ar[d]  \ar[r] ^{\hskip 1em p} &  \bH_g \ar@{=}[d]\\
& T_g\backslash \mathcal{T}_g \ar[d]  \ar[r]  &  \bH_g \ar[d] \\
&M_g= \Gamma_g\backslash \mathcal{T} _g  \ar@{^{(}->}[r]& \smpl{g}_\Z\backslash \bH_g.
} \end{xy}
\]
}
2. The homogeneous space $\smpl{n,\R} / \ugr{n}$ equals the Siegel upper half space $\bH_{n}$ parametrising  {marked\footnote{See the previous remark.}} polarized abelian varieties  of dimension $n$ with  polarization given by  the symplectic form $\mathbf J$. The action of the symplectic group  is  given by 
\[
Z\mapsto   \langle T \rangle \cdot Z:= (A+CZ)^{-1}(B+DZ),\quad T=\begin{pmatrix}A& B\\ C& D\end{pmatrix}\in  \smpl{n,\R} 
\]
Writing $Z=X+\ii Y$, one verifies that
\[
  \langle T\rangle \cdot  \ii  \mathbf {1}_n  = Z,\qquad T:=\begin{pmatrix} Y^{-\half} & Y^{-\half} X\\ 0& Y^\half\end{pmatrix}.
\]
Moreover, the corresponding complex structure and compatible metric are  given by
\[\begin{array}{lcl} J & =  &
T^{-1}\mathbf{J} T = \begin{pmatrix} -XY^{-1} & -Y-XY^{-1}X\\ 
Y^{-1} & Y^{-1}X\end{pmatrix}\\
 b_{\omega,J}  &=  &\Tr T T = 
\begin{pmatrix} Y & Y^{-1} X\\ 
XY^{-1} & XY^{-1} X+Y\end{pmatrix}.
\end{array} 
\]

\end{rmks}
\section{Examples}\label{sec:Examples}

\subsection{Cohomology of Riemannian Manifolds}
 \label{ssect:FirstEx}

\subsubsection*{Middle Odd Cohomology }\footnote{This example has been known a long time under the name of \emph{Lazzeri's jacobian} and has been studied in detail by Elena Rubei \cite{lazzeri}.} 

Let $(X,g)$ be a compact oriented Riemannian manifold of dimension $4q+2$ and consider the  lattice of integral middle cohomology modulo torsion, $H^{2q+1}(X)_\Z$ equipped with the cup-product pairing $\omega^-$. On $H^{2q+1}(X)_\R$ we put the Hodge metric (Def.~\ref{HodgeMetric}).
The abelian variety 
\[
J^{q+1}(X,g):=J(H^{2q+1}(X)_\Z,b^{(g)},\omega^-)
\]
 is an intrinsic invariant of the pair $(X,g)$. The polarization is a principal polarization.  The complex structure is just coming from the $*$ operator on middle cohomology. Note that the pair $(b^{(g)},\omega^-)$ is coherent in the sense of Def.~\ref{DefCoherent}.2.
 
 When $X$ is a compact topological  surface ($q=0$)  by Remark~\ref{OneDimCase}  this construction gives  back the classical jacobian of the associated Riemann surface and in fact, since the $*$ operator in cohomology is the   complex structure on the tangent space of the classical jacobian and the principal polarization comes from the usual cup product, by the classical Torelli theorem, giving $J^{1}(X,g)$ is completely equivalent to giving the isomorphism class of the associated Riemann surface.

\subsubsection*{Odd Cohomology} 

Suppose   that $(X,g)$ is an even $2d$-dimensional  compact connected oriented Riemannian manifold. On  $H^*(X)_\R$ we put the Hodge metric $b^{(g)}$ (Def.~\ref{HodgeMetric}) while on $H^-(X)$ we have the intersection pairing $\omega^-$. It is integral and unimodular on $H^-(X)_\Z$ and thus one has  unique complex structure associated to $H^-(X)_\R$. Indeed, this complex structure is given by the Hodge $*$ operator which on odd cohomology indeed is a complex structure. The pair $(b^{(g)},\omega^-)$ is a coherent pair  in the sense of Def.~\ref{DefCoherent}.2.

Note that   $\omega^-$   pairs $H^{2j+1}(X)_\Z$ perfectly to $H^{2d-2j-1}(X)_\Z$. 
So depending on whether $d$ is odd or even two cases arise both of which are coherent:
\begin{enumerate}
\item $\dim(X)=4q=2d$:
\[
J^-(X,g):=J(H^-(X)_\Z,b^{(g)},\omega^-)= J^1 \times J^2 \cdots\times  J^{q},\ 
\]
where $J^k=J^k(X,g)=J(H^{2k-1}(X)_\Z\oplus H^{4q-2k+1} (X)_\Z,b^{(g)},\omega^-)$.  
\item  $\dim(X)=4q+2$:
\[
J^-(X,g)=J(H^-(X)_\Z,b^{(g)},\omega^-)= J^1 \times J^2 \cdots\times  J^{q+1},
\] 
where all but the last factors are as before and the last factor is the invariant $J^{q+1}(X,g)$  {we previously considered}. \footnote{Note that except for the middle dimension, the odd cohomology groups are not in general even dimensional, hence the need to combine $H^{2k-1}$ with its ``dual'' $H^{4q-2k+3}$. This is in contrast with odd cohomology for K\"ahler manifolds. See Remarks~\ref{Confront} (2).}
\end{enumerate}

\subsubsection*{Even Cohomology}
Next,  consider the even cohomology.  Here we assume that $\dim X=4q+2$. Again we take the Hodge metric on $H^+(X)_\R$ but now we take the symplectic  form $\omega^{ +}$.

The resulting principally polarized abelian variety is denoted $J^{+}(X,g)$. It  is a product whose factors are  principally polarized abelian varieties associated to  two summands of the form $H^{2k}(X)\oplus H^{4q-2k+2}(X)$. The   complex structure on such a factor comes from
$\begin{pmatrix}  0 &  -*\\ * & 0\end{pmatrix}$ since $*^2=\id$ on even  cohomology. Again this gives a coherent example (Def.~\ref{DefCoherent}.2).

\subsubsection*{Twisted versions}

We now consider a twisted version of the above.   Given a vector space  isomorphism $ \gamma$ of  $H^*(X)$ preserving even and odd cohomology, put: 
  \begin{eqnarray}
 H^{\pm }(X) \times H^{\pm}(X)  \ni (x,y) \!\! & \mapsto &\!\!
 \omega^\pm_\gamma(x,y):=  \omega^\pm (\gamma(x), \gamma(y )  ). 
 \label{eqn:TwistProdOne}
 \end{eqnarray}
This  new pairing is again symplectic, but it  need no longer be integral on $H^{\pm}(X)_\Z$ but one may choose  $N=N(\gamma)$  to be a minimal  integer such that $N\omega^\pm_{\gamma}$ becomes integral.  This gives a canonical integral twist. 

We can likewise use an $\R$-vector space isomorphism  $\tau$ of $H^*(X)_\R$ to modify the Hodge metric: 
we define
 \begin{equation}\label{eqn:TwistProdTwo}
 b^{(g)}_{\tau}(x,y):= \int_X  \tau(  x ) \wedge  * \tau (y) 
 \end{equation}
Of course now in general the pairs  $(b^{(g)}_{\tau}, \omega^\pm_{\gamma})$ need no longer be coherent, but Lemma~\ref{Compatible} states: 
\begin{lemma} \label{WhenCohPairs?}  The pair $(b^{(g)}_{\tau}, \omega^\pm_{\gamma})$ is coherent if and only if $\gamma^{-1}\tau$ is  a positive multiple of an  $\omega^\pm$-symplectic map.
\end{lemma}
So we have: 
 \begin{prop}\label{AVforIndexLikesTwo} Let $(X,g)$ be a compact  {oriented} Riemannian manifold of dimension $2 \bmod 4$ and let   $ \gamma$ be a $\Q$-vector space isomorphism of  $H^*(X)$ preserving even and odd degree, and $\tau$ any $\R$-isomorphism of $H^*(X)_\R$. Using the notation of Definition~\ref{DEFN}, the above construction yields an  abelian variety
\[
J(H^{\pm}(X)_\Z, b^{(g)}_{\tau}, N(\gamma)\omega^\pm_\gamma),
\]
 canonically associated to $(X,g,\gamma,\tau)$.
 \\
 The pair $(b^{(g)}_{\tau}, \omega^\pm_\gamma)$ is coherent if and only if $\tau\gamma^{-1}$ is  a positive multiple of an  $\omega^\pm$-symplectic map. 
\end{prop}

\subsubsection*{Moduli Interpretation}   For simplicity we only consider manifolds whose dimension is $2\bmod 4$. One can give a similar statement for any even-dimensional manifold when one restricts only to odd cohomology. 
\begin{thm} \label{CohAsModuli} Let  $X$ be a  fixed smooth  {compact oriented} manifold of dimension $ {2}d=4q+2$. Let   $\mathsf{Conf}(X)$ be the space     {of} classes of  conformal metrics on $X$ and $\mathsf{Conf}(H^*(X)_\R \times \Lambda^2 [H(X)^\vee ]^0)$ the set of conformal equivalence classes of pairs (metrics, symplectic forms) on  $H^*(X)_\R \times   H^*(X)$.  The space of symplectic  forms on $H^*(X)  $ is denoted  $\Lambda^2 [H(X)^\vee ]^0 $.  Let $G$ be the group of $\R$-vector space isomorphisms of $H^*(X)_\R$ which preserve even and odd degree classes, $G_\Q$  the rationally  defined isomorphisms and  let  $G_\omega$ be the subgroup of isometries preserving $\omega=\omega^++\omega^-$.
Finally, let  $ g_1:= \sum_{j=1}^q b_{2j-1}+\half b_{2q+1}$,  $g_2:=\half\dim H^+(X)$.
 
The construction of Prop.~\ref{AlgProp} applied to $( b^{(g)} , \omega)$ defines $r$ in   the following commutative diagram
\[
\begin{xy}
\xymatrix{
\mathsf{Conf}(X)\times \set{\omega}  
  \ar[rrd]^{p}     \ar[rr]^(0.4){b^{(*)}\times \id}   & &  \mathsf{Conf}(H^*(X)_\R \times \Lambda^2 [H(X)^\vee ]^0)     \ar[d]^r  \\
 && {\mathbb H}_{g_1} \times \bH_{g_2}.}
 \end{xy}  
\]
The map $p$ factors   over the inclusion
\[
\underbrace{\prod_{j=1}^q \bH_{b_{2j-1}}\times \bH_{ {\half} b_{2q+1}}}_{\text{ \rm odd rank}} \times \underbrace{ \prod_{j=0}^q \bH_{b_{2j}}}_{\text{ \rm even rank}}  \into \bH_{g_1}\times \bH_{g_2}.
\]
The group  $G\times G_\Q$ acts naturally  on the source of the map $r$, the subgroups  $\sett{(\lambda g,g)}{\lambda\in \R^+, g\in G_\Q}$  and  $\R^+\cdot G_\omega\times \set{1} $ maps $(b^{(g)},\omega)$ to another coherent pair.  
\end{thm}
\begin{rmks}\label{rmksonconf}
1)  The left hand space is a  Teich\-m\"uller space for $X$ and  $p$ should  be viewed as a substitute  period map in the setting of compact smooth manifolds of  dimension $2\bmod 4$.
\\
2)  If, moreover $X$ is a symplectic manifold with symplectic form $\omega$, the space $\mathcal{J}^\omega(X)$ of almost complex structures on $X$ tamed by $\omega$  embeds in  $\mathsf{Conf}(X)$ by the map $\beta(J)=   g_{J,\omega}$ and this subspace
is a better substitute for the moduli space in this case. The period map restricts to it and in fact, we shall make use of this remark in the  K\"ahler  setting. So, instead of the above diagram we should use
\[
\begin{xy}\xymatrix{
\mathcal{J}^\omega(X)   \times \set{\omega}   \ar[r]^-(.3){\beta \times 1}  \ar[rd]^{p} &  \mathsf{Conf}(H(X))\times  \Lambda^2  [H(X) ^\vee] ^0)    \ar[d]^r \\
& {\mathbb H}_{g_1} \times \bH_{g_2}.}  \end{xy} 
\]
\end{rmks}

 \subsection{$K$-Groups of Compact Smooth Manifolds}  \label{ssect:SecondEx}

We can give a variant of the examples in \S~\ref{ssect:FirstEx} using $K(X)$.  Recall, \eqref{eqn:ChernIsIso} that the Chern character:
\[
\text{\rm ch}: K(X) \to H^+(X)
\]
becomes  an isomorphism after tensoring with $\Q$.  So
\[
\Lambda(X)= \ch{K(X)}\subset H^+(X) 
\]
is a \emph{lattice}, i.e., a $\Z$-module of rank $\dim_{\Q}H^+(X) $.   The intersection pairing $\omega^+$  {as} well as the twisted pairing  \eqref{eqn:TwistProdOne}  induces  $\Q$-bilinear pairings on $\Lambda$, i.e.,
\begin{equation}\label{eqn:LambdaMult}
 \omega^+(\ch \xi  ,\ch \eta )= \int_X \ch{\xi\otimes\overline{\eta}}.
\end{equation}
In the framework of  $K$-groups we twist by vector space isomorphisms $\gamma$ defined as  multiplication  by  a unit    $\mathbf{b}$ in the ring $H^{4*}(X)$.   {Recall that  $\iota$ is the involution on $H^{2*}$ which on $H^{4*}$ is the identity and on $H^{4*+2}$ minus the identity. So $\mathbf{b}$ is invariant under $\iota$  and    so  multiplication with it    commutes with $\iota$.   Hence}  the twisted pairing becomes
 \[
 \omega^+_{\mathbf{a}}(\ch \xi ,\ch \eta )=\int_X \mathbf{a} \wedge\ch{\xi \otimes\overline{\eta}},\quad \mathbf{a}=\mathbf{b}^2.
 \]
 It can now happen that for specific $\mathbf{a}$    such a twisted   pairing becomes an integral pairing:
 \begin{dfn} A \emph{multiplier} is an element $\mathbf{a}\in H^{4*}(X)$ such that the pairing 
 \[
  \omega^+_{\mathbf{a}} : \Lambda(X)\times\Lambda(X) \to \Q
 \]
 is  \emph{integral}.  If $\mathbf{a}_0=1$, such a multiplier is called \emph{normalized}. 
 \end{dfn}
Now an interesting phenomenon occurs which is based upon  Poincar\'e duality:
 \begin{lemma}[\protect{\cite[3.7]{ah3}}]{Suppose that $X$ is \emph{torsion free\footnote{This  means that the integral cohomology $H^*(X;\Z)$ has no torsion.}.} Normalized multipliers always exist. Moreover, if} $\mathbf{a}\in H^{4*}(X)$ is a normalized multiplier, the pairing  $\omega^+_{\mathbf{a}} $  is unimodular on $\Lambda(X)$.
 \end{lemma} 
\begin{corr}  \label{KCohPairs} Let $(X,g)$ be a {torsion free}
compact oriented  Riemannian manifold with  $\dim X\equiv 2\bmod 4$ and let $\mathbf{a}\in H^{4*}(X;\Q)$   be a normalized multiplier and let $\mathbf{a}=\lambda \mathbf{b}^2 \in H^*(X,\R)$ for some $\lambda>0$. Then (using the notation of Definition~\ref{DEFN}) the pair $(b^{(g)}_\mathbf{b}, \omega^+_{\mathbf a})$ is coherent and 
\[
J(\Lambda(X),b^{(g)}_\mathbf{b}, \omega^+_{\mathbf a})
\]
is a principally  polarized abelian variety.
\end{corr}
\begin{exmples}\label{worldex} 
(1)   If $X$ is {torsion free and} spin,  there is a \emph{canonical} choice for a principally polarized abelian variety associated to the spin structure.  {Indeed, by  the Index Theorem  of  Atiyah and Singer (see \ref{ATSIT}), in view of \eqref{eqn:LambdaMult} we have}:
 \[
  \omega_{\hat A(X)}^{  +}: \Lambda(X)\times \Lambda(X) \to \Z,\quad (\xi,\eta)\mapsto \int_X\hat A(X)\ch{\xi\otimes\overline{\eta}}=\text{index}(\dirac  D_{\xi\otimes\overline{\eta }}),
  \]
 and the result follows since $\hat A(X)$ is normalized: it starts off with $1\in H^0(X)$.
  In order to get a coherent pair (Hodge metric, twisted cup-pairing), by  Cor.~\ref{KCohPairs},   we may take for $\mathbf{b}$  any multiple of  $\sqrt{\hat A}$, such as \footnote{The square root is unique and belongs to $H^{4*}(X)$ since  $\hat A(X)$ starts with $1$.} 
  \[
\mathbf{b}= 2\pi \sqrt{\hat A}
\]
which is the Witten-Moore choice  from \cite{world}.
 
Note that one can take $\mathbf{a}=1$ if $c_1(X)=0$, for instance if $X$ is a complex torus, or, more generally, any Calabi-Yau manifold.
\\
 (2)   If $X$ is a complex    manifold with $c_1(X)=0$   there is  {a priori}  \emph{another}  canonical choice for a principally polarized abelian variety. 
 {One  takes $ \mathbf{a}=\text{\rm td}(X) $ and for $\mathbf{b}$ one takes any multiple of $\sqrt{\mathbf a}$. The Todd genus is known to take values in $H^{4*}(X)$ if and only if $c_1(X)=0$. Much more is true:  by the calculations in \cite[1.7]{hir} and especially formula (12) in loc. cit. the Todd and $\hat A$--genus coincide in this case. In particular this \emph{does not give a new example!}}
\end{exmples}

These examples can also be considered with moduli, e.g. varying metrics in examples (1) and (2) and varying the complex structure in example (3).  Let us formulate the final result in a  setting which is common to all examples:
\begin{thm} \label{MainThmForKGroups} Let  $(X,\omega)$ be a {torsion free} compact symplectic oriented manifold  of dimension $4q+2$. Let $\Lambda(X)=\ch{K(X)}\subset H^+(X)$.  The set of normalized multipliers  is denoted
$H^{4*}(X)^\times_{\rm norm}$ and $g:= \text{\rm rank }  \Lambda(X)$.
Recall    the embedding  $ \mathcal{J}^\omega(X)  \into \mathsf{Conf}(X)$  (cf. also Remark~\ref{rmksonconf} 2). Using it,  
we have a commutative diagram:
\[
\begin{xy}\xymatrix{
\mathcal{J}^\omega(X)     \times  H^{4*}(X)^\times_{\rm norm}    \ar[r]^(0.45){\beta} \ar[rd]^{p} &  \mathsf{Conf}(\Lambda(X)_\R)\times  \Lambda^2 [\Lambda(X)^\vee ]^0    \ar[d]^r \\
& {\mathbb H}_{g}.}  \end{xy}  
\]
Here $\beta(J,\mathbf{a})= (b^{g_{J,\omega}},\omega_{\mathbf{a}}^{  +})$.
\end{thm}

\subsection{Hodge Structures} \label{sec:HS}

\subsubsection*{Odd weight} In this situation the underlying rational vector space has even dimension, say $\dim W=2g$. Above we saw that the Weil operator $C_W$  is a complex structure and so  it then defines a  weight $1$ Hodge structure, say $V$. If $W$ is polarized by $Q$ it is clear that $V$ is polarized by $Q$ (being a polarization only depends on the Hodge structure through the Weil operator). If $(W,Q)$ is an integral polarized Hodge structure, the corresponding polarized abelian variety is denoted $J(W,Q)$. It is the so-called \emph{Weil jacobian} \cite{weiljac}. 
\begin{rmk}
The complex structures $\pm C_W$ on $W$ are characterized by the $\pm \ii$-eigenspaces being the direct sums of the Hodge $W^{p,q}$-spaces with $ p\equiv k \pmod 2 $ and $ p\equiv k +1 \pmod 2$ respectively, i.e., one places the  Hodge spaces  alternatingly in the two different eigen-spaces.
There is an obvious different choice for the complex structures  $\pm \tilde C_W$   in which the two $\pm \ii$- eigenspaces  are given by the sum of the $W^{p,q}$ with $p\ge k $ and $p<k$ respectively: the first (second) half of Hodge spaces form the first (second) eigenspace. This  last choice gives the Griffiths intermediate jacobian $\tilde J(W,Q)$.
For this choice the compatible metric (replacing $C_W$ by $\tilde C_W$) is no longer positive definite but in general indefinite. It no longer gives a Riemann form on the torus, so $\tilde J(W,Q)$ need not be an abelian variety.
\end{rmk}

 In terms of homogeneous spaces, the Griffiths domain    for the type $W$-structures is 
$\smpl{g}/H$, where $g=h^{2\ell -1,0}+\cdots+ h^{\ell,\ell-1}$ and
\[
H=\ugr{h^{2\ell -1,0}}\times\cdots\times \ugr{h^{\ell,\ell-1}},\qquad k=2\ell-1
\]
while the one for the $V$-type structures is $\smpl{g}/\ugr{g}$. The assocation $W\mapsto V$ is induced by the natural map
\begin{equation}\label{eqn:firstBarPsi}
\bar\psi : \smpl{g}/H \to \smpl{g}/\ugr{g}\simeq \germ h_g
\end{equation}
and is well known to be in general neither holomorphic nor anti-holomorphic \cite[3.21]{periodarticle}. 
This is easily illustrated in weight $3$ as follows. The Griffiths domain para\-metrizes Hodge flags 
$F^3\subset F^2$ inside $W_\C$ satisfying the two Riemann conditions. The subspace $F^3+F^1 \cap \overline{F}^2 $ then is a $g=(\half\dim W$)-dimensional  isotropic subspace of $W_\C=V_\C$, i.e., satisfies the first Riemann condition (Lemma~\ref{frc}) and it satisfies also the second Riemann condition and we have $\bar\psi (F^3,F^2)=  F^3+F^1 \cap \overline{F}^2 $ which obviously is non-holomorphic (consider the Pl\"ucker coordinates).

Even  in geometric situations when there is a holomorphic period map $p:M \to \smpl{g}/H $, the composition $\bar\psi\comp p$ is seldom holomorphic as shown Griffiths' calculation \cite[Proposition 3.8,  and I, 1.3 (b)]{periodarticle}  where  the triple product of the universal family of  elliptic curves over the upper half plane is considered in detail. Here  $M= \germ h_1^3$ and the composed map becomes  $\bar\psi\comp p: \germ h_1^3 \to \germ h_3$.

In the case of simple polarized weight $3$ Hodge structures $W$ with $h^{3,0}=1$ it is known that if $W$ is CM, i.e., its Mumford--Tate group is a torus, then also the associated Weil Intermediate jacobian has CM \cite{borcea}. This amounts to saying that $\bar\psi$ preserves complex multiplication. More generally, $\bar\psi$ preserves additional endomorphisms of the Hodge structure $W$. 
In fact, this is obvious from a Tannakian point of view, since the Mumford--Tate group of $W$ resp. $V$ is the Tannaka group of the rigid tensor subcategory generated  
by $W$ resp. $V$ and its tensor powers. 
A similar remark applies in the following section. 

\subsubsection*{Even Weights}

Here we take $V=W\oplus W^\vee(-k)$. It has a natural complex structure given by the linear map
\[
J(x+ \hat Q y)= \hat Q C_W(x) -C_W(y)
 \]
 and so, by the above, defines a weight one Hodge structure on $V$ with $C_V=J$. The polarization $Q$  defines a polarization $q$ as follows:
 \[
 q(x_1+\hat Q y_1, x_2+\hat Q y_2)= -Q(x_1,y_2)+Q(y_1,x_2).
 \]
 By construction $q$ is skew-symmetric and it is a standard verification that $q$ is skew and satisfies the Riemann bilinear conditions.
 For instance,  Lemma~\ref{fbl}  shows that the first bilinear relation can be tested by showing  that $q$ is $C_V$-orthogonal which is the case since
 \[
 \begin{array}{lcl}
 q(C_V(x+\hat Q y),C_V(x'+\hat Q y')) &=& q(\hat Q C_Wx-C_Wy,\hat Q C_W x'-C_Wy')\\
 &=& Q(C_Wx,-C_Wy')+Q(C_W {y},C_W x')\\
&=& -Q(x,y')+Q(y,x')\\
 &=&q(x+\hat Q y,x'+\hat Q y').
 \end{array}
 \]
If   $(W,Q)$ is an integral polarized Hodge structure, $(V,q)$ is integrally polarized and    if $Q$ happens to be unimodular, also $q$ is unimodular. The corresponding abelian variety $J(V,q)$ will also be denoted $J(W,Q)$.

To interprete the construction in terms  of  Griffiths domains we need to introduce the homomorphism
\begin{eqnarray}
\psi  :  \so{W,Q}& \to &\smpl{W\oplus W^\vee(-k),q} \label{eqn:psi} \\
  f & \mapsto& \psi(f),\qquad \psi(f)(x+\hat Q y)= f(x)+\hat Q (f(y)).\nonumber
\end{eqnarray}
To see this is well defined we need to verify that $\Phi:=\psi(f)$ is indeed symplectic:
\begin{eqnarray*}
q( \Phi(x+\hat Q y), \Phi (x'+\hat Q y'))& = & q(f(x)+\hat Q f(y), f(x')+ \hat Q f(y'))\\
&=& -Q(f(x),f(y'))+Q(f(y),f(x'))\\
&=& -Q(x,y')+Q(y,x')= q(x+\hat Q y,x'+\hat Q y').
\end{eqnarray*}
Now write $k=2\ell$. Let $p=\sum _j h^{2\ell-2j ,2j}$ and $q= k- p  $. Then $\so{W,Q}\simeq \so{p,q}$ and
\[
H=  \ugr{2h^{2\ell,0}} \times\cdots \ugr{2h^{\ell+1,\ell-1} } \times \so{h^{1,1}}
\]
\begin{lemma} \label{NewPeriods} Let $\so{W,Q}/H$ be the Griffiths domain for polarized Hodge structures of type $(W,Q)$.  Suppose $\dim W=g$.
Then the map $\psi$ from \eqref{eqn:psi} induces a diagram
$$
\xymatrix{ 
\so{W,Q}/H  \ar^{\pi}[dr]   \ar[rr]^{\bar\psi} &  & \smpl{g}/\ugr{g}  \\
& \so{W,Q}/K  \ar[ur]  , & }
$$
where $K\subset  \so{W,Q}$ is the unique maximal compact subgroup containing $H$, $\pi$ the natural map and $\tilde \psi$ the induced map.
\end{lemma}
\proof  Let $F \in \so{W,Q}/H $   correspond to the given polarized Hodge structure $(W,Q)$.  Let 
\[
h: \mathbf{G}(\R) \to \ogr{W,Q}
\]
be the representation which gives this Hodge structure. Then $H$ is the commutant of the  Mumford-Tate group of the given Hodge structure (the Zariski-closure of the image of $h$). Indeed, $H=\sett{g\in  \ogr{W,Q}}{  g h(z)= h(z) g  \text{  \rm for all } z\in \mathbf{G}(\R)}$.  Now, by construction $\psi\comp C_W= C_V\comp \psi$ and hence $\psi$ sends the commutant of $C_W$ to $\ugr{g}$, the commutant of $C_V$.   But   $H$, the commutant of the  Mumford Tate group of $W$ is   contained in  the commutant of $C_W$   and so $\psi$ sends  $H$ to  $\ugr{g}$.

Finally we have to show that  $\bar\psi$ factors over the quotient of $\so{W,Q}$ by the maximal compact subgroup $K$. To see this consider $\psi(K)+\ugr{g}\subset \smpl{g}$. It is a compact Lie subgroup containing $\ugr{g}$ and since the latter is already maximally compact, $\psi(K)\subset \ugr{g}$ and hence $\bar\psi$ factors over the quotient map $\pi: \so{W,Q}/H  \to \so{W,Q}/K$.
\qed\endproof

As in the odd weight case the map $\bar \psi$ is in general not holomorphic. We give a proof for  weight $2$. Abbreviate  $F=H^{2,0}$. Then $H^{1,1}= (F+\bar F)^\perp$ (with respect to the polarization). The Griffiths domain in this case is an open subset of the Grassmann variety of $\dim S$-dimensional linear subspaces of $W_\C$. Define the two maps
\begin{eqnarray*}
\varphi_{\pm}: V_\C  & \to & W_\C+W_\C^*(-2)=V_\C\\
 z &\mapsto &  z \pm  \ii \hat Q (z). 
\end{eqnarray*}
The images are the $\pm \ii$-eigenspaces for $C_V$ and hence the map $\bar \psi$ comes from the map
\[
F\mapsto \varphi_+(F+\bar F) +\varphi_-(F+\bar F)^\perp \subset V_\C.
\]
Considering Pl\"ucker coordinates, one sees that this map is   neither  holomorphic nor anti-holomorphic as soon as $h^{1,1}\not=0$.

\begin{rmq} In geometric situations the composition of the period map with $\bar\psi$ is in general non-constant (even not holomorphic and not anti-holomorphic).
\end{rmq}

We illustrate this  with  the following 
\begin{exmple}  Let $E=E_\tau$ be the elliptic curve $\C/\Z+\Z\tau$ and let $\alpha,\beta$ be the two cycles coming from the two lattice generators $\set{1,\tau}$. Then $E\times E$ has a natural principal polarization inducing one on $W=H^1(E)\otimes H^1(E)$. The cycles $\alpha\otimes\alpha$, $\alpha\otimes\beta$, $\beta\otimes\alpha$, $\beta\otimes\beta$ define a lattice $W_\Z\subset  W$ on which the polarization $Q$ is integral (even unimodular).
Let  $\set{(\alpha\otimes\alpha)^*, (\alpha\otimes\beta)^*, (\beta\otimes\alpha)^*,(\beta\otimes\beta)^*}$ be the dual basis. Then 
\begin{equation}\label{eqn:ActionHatQ}
 \alpha\otimes\alpha  \overset{\hat Q}{\iff} (\beta\otimes\beta)^*,\qquad  \beta\otimes\alpha  \overset{\hat Q}{\iff} (\alpha\otimes\beta)^*.
\end{equation}
If $\omega=dz$ is the normalized $1$-form on $E$ with periods $1$ and $\tau$, the    Hodge structure on $W$ has $h^{2,0}=1$ with basis $\omega\otimes\omega$ and period matrix  
\[
 (\tau^2,\tau,\tau,1) \in   \C^4= W_\C
\]
spanning the line $F\subset W_\C$. The weight one Hodge structure on $V$ then is given by calculating the periods of $\omega\otimes\omega$ with respect to a suitable basis for $ \varphi_+(F+\bar F) +\varphi_-(F+\bar F)^\perp$. 

As a basis for $V$ we take 
\[\set{ \alpha\otimes\alpha, \alpha\otimes\beta, \beta\otimes\alpha,\beta\otimes\beta,
 (\alpha\otimes\alpha)^*, (\alpha\otimes\beta)^*, (\beta\otimes\alpha)^*,(\beta\otimes\beta)^*}
 \]
 Then $F+\bar F \subset W_\C$ is given by the matrix
 \[
M:= \begin{pmatrix}
 \tau^2 & \tau & \tau & 1\\
 \bar\tau^2 &\bar\tau & \bar\tau & 1
\end{pmatrix}.
 \]
Define an involution $\iota$ on $(4\times 2)$-matrices: exchange column 1 and 4 as well as column 2 and 3. In view of \eqref{eqn:ActionHatQ}, the subspace $\varphi_+(F)\subset V_\C$ is given by the matrix $(M, \ii\iota(M))$. 

The matrix for $Q$ in the given basis then is the $(4\times 4)$ anti-diagonal  matrix with $1$ on the antidiagonal and hence  $(F+\bar F)^\perp$ is  given by calculating  a basis for $\ker M$ and then applying $\iota$. We find
\[
N:=\begin{pmatrix}
|\tau|^2 & -(\tau+\bar\tau) & 0 & 1\\
0 & 1 & -1 & 0
\end{pmatrix}.
\]
Then $ \varphi_+(F+\bar F) +\varphi_-(F+\bar F)^\perp$ is given by the block matrix
\begin{equation}  \label{eqn:Block}
B:=\begin{pmatrix}
M & \ii \iota(M) \\
N & -\ii \iota(N)
\end{pmatrix}.
\end{equation}
One calculates $\det \begin{pmatrix}
M   \\
N  
\end{pmatrix}=(\tau-\bar\tau)(\tau^2+6|\tau|^2+ \bar\tau^2)$ and similiarly, one finds among  the other non-zero Pl\"ucker coordinates $-\ii (\tau-\bar\tau)(\tau+\bar\tau)^2$ (e.g., in the previous determinant, replace the first column by the last column of the block matrix \eqref{eqn:Block}). Since the quotient of these two equals $-\ii \displaystyle\frac {\tau^2+6|\tau|^2+\bar\tau^2}{(\tau+\bar\tau)^2}$ the period map composed with $\bar\phi$ is a \emph{non-constant} map which is neither holomorphic nor anti-holomorphic.

\end{exmple}

\subsection{Cohomology of  K\"ahler manifolds} \label{ssec:CohKaehler}
Let $(X,\omega)$ be a projective manifold of dimension $d$ with integral K\"ahler class $[\omega]$.   We have seen (Example~\ref{CohKaehler}) that the primitive cohomology groups $H^k(X)_{\rm prim}$ as well as the full cohomology groups $H^k(X)$ have a natural weight $k$ Hodge structure polarized by  the form $Q_\omega$ defined by  \eqref{eqn:PolarAll}. We can then apply the two constructions  of \S~\ref{sec:HS}. In fact, for primitive $k$-cohomology  these match exactly  the constructions  we described in the general setting of compact smooth manifolds of dimension $2 \bmod 4$ as given in \S~\ref{ssect:FirstEx}.

\subsubsection{Odd Cohomology}  For odd rank $2k+1$ we have the  Weil jacobian 
\[
J^{k+1}(X)=J(H^{2k+1}(X)_{\rm prim}, Q_\omega).
\]
In this case  $(b^{g_\omega}, Q_\omega)$  is a coherent pair  with  the complex structure given by   the  Weil operator  $C$.  This also works on $H^{ 2k+1}(X)$ except that we now have to multiply $Q_\omega$ with a certain  integer  $N_{k,\omega}$ yielding    a polarized abelian variety
\[
J(H^{2k+1}(X)_\Z, b^{g_\omega}, N_{k,\omega} Q_\omega).
\]
This torus is isogenous\footnote{Two tori $A$ and $B$ are said to be \emph{isogenous} if there is a surjective group homomorphism $A\onto B$ with finite kernel. Despite the apparent asymmetry  in the definition this does define an equivalence relation.} to a product of tori coming from the primitive pieces, i.e.,  with $\sim$ denoting isogeny, we have
\[
J(H^{2k+1}(X)_\Z, g_\omega, N_{k,\omega} Q_\omega)
\sim \prod_{\ell=1}^ {k+1}   J^\ell (X).
\]

\subsubsection{Even Cohomology} Here we have to  restrict to odd dimensional complex varieties in order that the real dimension $2d$ be $2\bmod 4$. The even cohomology can now be given as a direct sum $W\oplus W^\vee$ where $W$ is the sum of the first half of the even  cohomology groups. Then, using Poincar\'e duality  $W^\vee$ is the sum of the last half of the even cohomology groups. Incorporating the Hodge structure one should pair $H^{2k}(X)$ with $H^{2d-2k}(X)(d-2k)$ and the twisting cup-product pairing should be used for the symplectic form.
The construction we have given for abstract Hodge structures in \S~\ref{sec:HS} for even weight, then is the same as the third example  from  \S~\ref{ssect:FirstEx} (even cohomology).

\subsubsection{Moduli}
Let $(X,J,\omega)$ be almost K\"ahler. This means  that $\omega$ is a non-degenerate real two form, $J$ an almost complex structure so that $g_J+\ii\omega$ is a hermitian metric and $d \omega=0$. We fix $\omega$ and let the almost complex structure $J$ vary. The coherent pairs $(g_J,\omega)$ up to conformal equivalence form the space $\mathsf{Conf}^\omega$. Restricting to \emph{integrable} complex structures we get $\mathsf{Conf}^\omega_{\rm int}$, corresponding to the true  K\"ahler metrics. The fixed class $[\omega]\in H^2(X)_\R$ can be used to define primitive cohomology and hence we have a polarized Hodge structure on $H^*_{\rm prim}(X)_\R$ parametrized by a   Griffths period domain $D$ and period map $p$. 
We arrive at the following moduli interpretation:

\begin{thm} \label{CohKaeherAndModuli} Let $(X,\omega)$ be a compact K\"ahler manifold of odd complex dimension.  Let $g_1:=\half \dim H^-(X)_{\rm prim}$ and  $g_2:=  \dim H^+(X)_{\rm prim}$.
 With the map $\bar \psi$  the one from \eqref{eqn:firstBarPsi} (odd cohomology) combined with the one from   Lemma~\ref{NewPeriods} (even cohomology) we have the following commutative  diagram
\[
\begin{xy}\xymatrix{
\mathsf{Conf}^\omega    \ar[d]^{} & \ar@{_{(}->}[l] \mathsf{Conf}^\omega_{\rm int}  \ar[r]^{} \ar[d]^{} \ar[rd]^{} & 
\mathsf{Conf}(H^*_{\rm prim}(X)_\R)  \ar[d] \\
\mathcal{J}^\omega(X) & \ar@{_{(}->}[l]  \mathcal{J}^\omega_{\rm int}(X) \ar[d]^{} \ar[r]^{p}  & {\mathbb H}_{g_1}\times \bH_{g_2}. \\
& D \ar[ur]^{\bar \psi} & 
} \end{xy}
\]
Moreover, if $\dim_\R X=4q+2$ and $q_j:= \dim H^j(X)_{\rm prim}$, the map $p$ factors over the inclusion
\[
\prod_{j=1}^{q+1} \bH_{\half q_{2j-1}}\times \prod _{j=0}^q \bH_{q_{2j}}\into \bH_{g_1}\times \bH_{g_2}.
\]
\end{thm}
In this theorem, the space $ \mathcal{J}^\omega_{\rm int}$ of complex structures  should be considered as  a sort of Teichm\"uller moduli space, and the map $p$ is a sort of period map which in this case  happens to factor over the Griffiths domain $D$.

\begin{rmks} \label{Confront} (1) A similar statement holds for all of cohomology but it does not add new information since the new tori are products of the ones gotten from primitive cohomology \\
(2) Note that  although the  diagram of Theorem~\ref{CohKaeherAndModuli} is very similar to the diagram in Theorem~\ref{CohAsModuli} which applies to the cohomology of \emph{any} compact Riemannian manifold of dimension $d=4q+2$, there is one crucial difference: the tori occurring in the latter situation need to be constructed by combining $H^{2k+1}(X)$ and its ``dual'' $H^{4q-2k+1}(X)$ since $\dim H^{2k+1}(X)$ need not be even dimensional. In the K\"ahler setting this is true however and the tori split in sub-tori coming from the various primitive pieces. To be precise, with $\sim$ denoting isogeny, we have
\[
J(H^{2k+1} (X)_\Z\oplus H^{4q-2k+1} (X)_\Z,b^{(g_\omega)},\omega^-) \sim \left[ \prod_{\ell=1}^ {k+1} J^\ell(X)\right]\times \left[ \prod_{\ell=1}^ {k+1} J^\ell(X)\right].
\]
\end{rmks}

\section{Special Theta Functions and Ramond-Ramond Fields} \label{MainResult}
 \subsection{Some reminders}
Let $J=V/\Lambda$ be an abelian variety with principal polarization given by the unimodular integral and positive $(1,1)$-form $\omega$. 
In this section let $g=\dim_\C V$ be the complex dimension of $J$. The set of line bundles $L$ on $J$  with $c_1(L)=\omega$ is a principal space under the Picard torus of $J$ which is isomorphic to $J$.
To single out a line bundle having $\omega$ as first Chern class one traditionally uses multipliers  for $\omega$:
\begin{dfn}  \label{dfn:mult} A \emph{multiplier for} $\omega$ is a  function 
$\alpha: \Lambda \to   \mathrm{U}(1)$
for which 
\begin{equation}\label{eqn:mult}\alpha(x+y)=(-1)^{\omega(x,y)}\cdot \alpha(x)\alpha(y).
\end{equation}
\end{dfn}
A multiplier is entirely specified by its values on a symplectic basis and any value in $ \mathrm{U}(1)\simeq S^1$ can be taken so that the above set is  a topological torus of dimension  $2g $ as should be the case. 
Indeed,  line bundles with given polarization $\omega$ are in 1-1 correspondence with multipliers for $\omega$. See for example \cite[Chap I.2]{mum}.
A choice of a symplectic basis  $\mathsf{B}:=\set{e_1,\dots,e_g, f_1,\dots,f_g}$  for  $\Lambda$ with respect to $\omega$ singles out a specific line bundle:  the one for which $\alpha(b)=1$ for all $b\in \mathsf B$. In particular, the corresponding multiplier  takes its values in the
 subgroup $\set{\pm 1$} of  U$(1)$. There are exactly $2^{2g}$ such  {``special''} line bundles since one 
may choose $\alpha(b)\in \set{\pm 1}$ for   every individual  $b\in \mathsf{B}$ separately  {; these 
correspond classically to {\em theta functions with characteristics}.}

In a more explicit fashion,  take a holomorphic basis for $V$, or, equivalently, a  basis  for  the space of holomorphic $1$-forms on $J$ chosen in such a way that the rows in the matrix
\[
\Omega=( \mathbf{1}_g, Z),\quad \Tr Z=Z, \,\im (Z)>0 
\]
are the periods of this basis  with respect to $\mathsf{B}$. Choose $\theta \in \Lambda $ such that $\alpha(y)=(-1)^{\omega(\theta,y)}$ for all $y  \in \Lambda $ which is possible since we have a symplectic basis. Then, setting 
\begin{eqnarray*}
\Lambda_1=\bigoplus \Z e_j && \Lambda_2:= \bigoplus \Z  f_j \\
\theta= \theta_1+\theta_2,&&  \theta_i\in \Lambda_i\\
u= \half \theta_1 \bmod   \Lambda_1 \in\half  \Lambda_1/\Lambda_1,&&  v= \half \theta_2 \bmod \Lambda_2 \in \half \Lambda_2/ \Lambda_2 
\end{eqnarray*}
define 
\[
\Theta\left[ {u\atop v} \right] (z) := \sum_{x\in \Lambda_1+   u}\ \exp [\ii \pi \langle x, Z x \rangle]\cdot   \exp[2\pi \ii \langle x,z+  v\rangle].
\]
It is the classical theta function with  \emph{theta characteristic}  $(u,v)$. Here $\langle x, y\rangle=\Tr 
x \cdot y$ is the usual euclidean inner product on $\C^g$.  It is the non-zero holomorphic section,  unique  
up to a multiplicative constant   for  the (unique)  holomorphic line bundle on $J$ with such a special 
multiplier $\alpha$. Note that the classical theta function corresponds to $\alpha=1$, \emph{but this is not 
the one suitable for physics,{according to \cite[\S~3.1]{world}}, as we shall see in the next 
subsection}. 

\subsection{Ramond-Ramond fields} 

Continue with the example~\ref{worldex} constructed from a {torsion free} compact spin manifold  $(X,g)$. So on $\Lambda=\Lambda(X) \simeq K(X)$   there is a natural unimodular  symplectic  form $\omega$ given by  
\[\omega(x,y)= \omega^+_{ \hat A(X)} ( \ch x ,\ch y ).
\] 
Next,  one needs to assume that $\dim(X) \equiv 2 \bmod 8$.     By   Prop.~\ref{RealIndexThm}    one   then   has a homomorphism
$ j: \mathrm{KO}(X)  \to \Z/2\Z .$ If $x\in \mathrm{K}(X)$ is a virtual complex bundle $x\otimes\bar x $ is naturally an element of $\mathrm{KO}(X) $ and so we get a homomorphism
\[
\alpha: K(X) \to \set{\pm 1},\quad x \mapsto (-1)^{j(x\otimes\bar x) }.
\]
One can   show  that it satisfies the required transformation law  \eqref{eqn:mult} to make it a multiplier for the form $\omega$.  
We  can now formulate the main result of \cite[\S 3]{world} in mathematical terms:
 \begin{prop}[\protect{\cite[\S~3.1]{world}}] For a {torsion free} compact spin manifold $(X,g)$ of dimension $2 \bmod 8$ consider the principally polarized abelian variety
  \[
 J(\Lambda(X),b^{(g)}_{2\pi \sqrt{\mathbf{a}}}, \omega^+_{\mathbf{a}}),
 \]
 where $\mathbf{a}=\hat A(X)$.  The map   $\alpha(x)= (-1)^{j(x\otimes\bar x)}$ is a multiplier  for $\omega$ and hence defines a  unique line bundle with first Chern class $\omega$ and multiplier $\alpha$. Let $\Theta\left[ {u\atop v} \right] (z)$ be the corresponding normalized theta function. Then the partition function for type II-A Ramond-Ramond fields on $X$ (see \S~\ref{PhysBack}) is given by $\Theta\left[ {u\atop v} \right] (0)/\Delta$ where $\Delta$ is the determinant \eqref{eqn:RegDet} for the non-zero modes on $X$.
    \end{prop}
In loc. cit. Moore and Witten argue that only this choice of the theta function gives the correct partition function. In particular, their result shows that   the partition function for the Ramond-Ramond fields on, say a 10--dimensional compact space-time  with spin structure    can be calculated completely  from a specific twist of the Riemann theta function which  is canonically associated to the spin structure.
 \appendix
\setcounter{section}{0}
  
\section{Clifford algebras}

Let $k$ be a field of characteristic $\not=2$ and let $V=(V,q)$ be a (finite dimensional) $k$-inner product space. Its tensor algebra  $\mathsf{T}V$ gets a natural $k$-inner product, also denoted by $q$. The unit $1\in k$ serves  as a unit in $\mathsf{T}V$.
Its \emph{Clifford algebra} is the following quotient algebra of dimension $2^n$ where $n=\dim V$:
\[ 
\cl V=\cl {V,q} := \mathsf{T} V /  \text{\rm ideal generated by }  \sett{x\otimes x+ q(x,x)\cdot 1 }{x\in V} . 
\]                                          
There is a natural map $c :V\to \cl V $, $ x\mapsto $class of $x$. The induced action $c(x):\cl V\to \cl V$ is  called  the \emph{Clifford action}.   
One easily shows that the pair $(\cl V,c)$ satisfies the following property: 
It is the unique pair $(C,c)$ consisting of a $k$-algebra  $C$ with unit together with  a  $k$-linear map $c: V \to C$ such that\footnote{The multiplication in $\cl V$ is written with dots.} 
\begin{equation}\label{eqn:ClifEq}
c (x)\cdot c(y)+c(y)\cdot c(x)= -2q(x,y)\cdot 1
\end{equation}
which is universal with respect to this property. 
\begin{propdef} A \emph{$\cl V$-Clifford module} $A$ is a $k$-algebra equipped with a \emph{Clifford action} of $V$, i.e., a linear map $c:V \to A$ satisfying \eqref{eqn:ClifEq}.
If the Clifford-action is $q$-skew-adjoint, one  says  that $A$ is \emph{self-adjoint}:
\[
q(v\cdot x,y)+q(x, v\cdot y)=0, \quad \text{\rm for all } v\in V,\, x,y\in A.
\]
\end{propdef}
The Clifford algebra is a twisted version of the exterior algebra:

\begin{lemma}  \label{WedgeIsCliff} $\cl V$, as a vector space  is isomorphic to $\Lambda V$, the exterior algebra. The Clifford-action on $\Lambda V$ is given by 
\begin{equation*}
c(x) \alpha = x\wedge \alpha - \iota(x)  \alpha 
\end{equation*}
where the linear map $ \iota(x)$ is the contraction with $x$. This makes $\Lambda V$ into a self-adjoint Clifford module.
 The bigrading given by odd and even degree in the  exterior product descends:
 \[
 \clif +V:=c\left(\Lambda^{+} V\right) ,\quad  \clif - V := c\left(\Lambda^{-} V\right).
 \]
 \end{lemma}
 \proof 
 Since $ \iota(x)$  is the  $q$-adjoint of the  map $  \alpha \mapsto  x\wedge \alpha $ clearly  $\Lambda V$ is  a self-adjoint Clifford module. 
 To see that one  gets an isomorphism, note that   
\begin{equation} 
\sigma: \cl V \mapright{\sim}  \Lambda V,\quad a\mapsto c(a)\cdot \mathbf{1} 
\end{equation}
  is a bijective $k$-linear map whose  inverse 
\[
c: \Lambda V \to \cl V
\]
can be  explicitly given as follows: let $\set{e_1,\dots,e_n}$ be an orthogonal basis for $V$, then send $e_{i_1}\wedge\cdots\wedge e_{i_k}$ to the element $e_{i_1}\cdots e_{i_k}$.   
 \qed\endproof
 
 From now assume that $\dim V$ is even. Clifford modules  turn out to be 
representations of  the \emph{spin group}  which in this case can be defined as  
\[
\spin V:=\sett{x_1\cdot \cdots \cdot x_{2k}\in \clif + V}{\|x_j\|=1, j=1,\dots,k}.
\]
To describe the basic irreducible Clifford modules extend $q$   bilinearly to $V_\C=V\otimes \C$. Then there exist maximal isotropic subspaces $H\subset V_\C$  with $\dim_\C H= \half \dim_\R V$.  If $V$ is oriented with orthonormal oriented basis $\set{e_1,\dots,e_{2m}}$ we can take for $H$ the subspace spanned by $e_{2k-1}+\ii e_{2k}$, $k=1,\dots,m$. Such $H$ is  an  \emph{oriented} maximal isotropic subspace. 
 Next,  introduce the \emph{spinor spaces}
\begin{equation}\label{eqn:SSpace}
\mathsf{S} =\sps{} V:= \Lambda H, \quad \sps  + V= \Lambda^+H ,\quad  \sps - V=\Lambda^-H. 
\end{equation}
The first, $\mathsf{S} $,  is clearly a complex Clifford module through the usual Clifford action given by Lemma~\ref{WedgeIsCliff}. It turns out to be an irreducible complex spinor representation.  On the other hand, the two spinor spaces  $\sps \pm V$  can be shown to be  irreducible as\emph{ real} representations of the spinor group. They are called the  \emph{half spinor representations}.

The metric  on the spinor space coming from the metric on $H$  induced by the hermitian form $(x,y)\mapsto q(x,\bar y)$   makes $\sps + V$ and $\sps - V$ orthogonal to each other and $\sps{} V$ is a self-adjoint Clifford-module.
One has:
\begin{prop} \label{TwistingSp} Every   complex $\cl V$-module $E$ is of the form $\sps{} V\otimes W$ where the \emph{twisting space} $W=\Hom_{\cl V}(\sps{} V,E)$ is a  complex vector space with trivial $\cl V$-action.  
\end{prop}

The spinor space is a $\Z_2$-graded complex $\cl V$-module. This is also the case for general Clifford  modules, but here one has to consider how the  \emph{chirality operator}  
\[
\gamma:= \ii^m e_1\cdots e_{2m} \in \cl V.
\]
acts:
\[
E^\pm:= \sett{e\in E}{  \gamma\cdot v= \pm v}.
\]
This is compatible with  action of $ \gamma $  on $\sps{} V$ since it turns out that  $ \gamma = \pm \textbf{1}$ on $ \sps \pm V$.   In particular, one has a $\Z_2$-graded action of $\cl V$ on $E$.

  \section{K-theory of Real Vector Bundles}
A reference for  this appendix is  \cite{at}, also contained as an appendix in \cite{at0}. \\
The Grothendieck group of real vector bundles on a manifold $X$ is denoted by $\text{\rm KO}(X)$. As in the complex case one sets $\text{\rm KO}^{-n}(X)=\text{\rm KO}(S^nX)$ and now there is periodicity of order $8$.  
 
There is   still another K-group   defined   for pairs $(X,\iota)$ where $X$ is a manifold and $\iota$ is an involution. One defines $\text{\rm KR}(X)$ as the K-group for complex bundles $E$ on $X$ admitting involutions covering $\iota$ and which are $\C$ anti-linear on the fibres. Again, there is a Bott-periodicity result, namely $\text{\rm KR}^*(X)\simeq \text{\rm KR}^{*+8}(X)$. The standard example is the total space of a \emph{real} vector bundle $E$ with involution $\iota$ given by $\id$ on the fibres.  This gives back $\text{\rm}KO(X)$. Another example is the Thom space $(BV,SV)$ of a Riemannian vector bundle $V$. Here the involution is the antipodal map. This space  figures  in  a very general form of the Thom isomorphism theorem which can be deduced from \cite[Theorem 6.2]{at2}.  We explain the latter theorem   in a simplified situation. Let $X$ be a compact differentiable manifold, $G$ a compact Lie-group acting trivially on $X$ and suppose we have a group homomorphism $\rho:G\to \text{\rm Spin}^c (8r)$. Moreover let $V$ be a  vector bundle on $X$ of rank $8r$ with spin$^c$-structure. Put a  $G$-module structure on $V$ through $\rho$. Then there is a natural isomorphism
\[
\varphi: \mathrm{KR}(X) \to \mathrm{KR}((BV,SV) \times_G  X).
\]
Specialize this to the case where $X$ is a spin manifold of dimension $(8r-m)$ so that $TX$ gets a spin$^c$-structure, let $G=\text{\rm Spin}^c (8r-m)$
and let $\rho:\text{\rm Spin}^c (8r-m) \into \text{\rm Spin}^c (8r)$ the embedding. Put $V=TX \oplus \R^m$. Then $(BV,SV) \times_G  X= (BX,SX)\oplus (B^m,S^m)$ where the involution on the second summand  is not the identity \emph{but the antipodal map}. Applying the periodicity \cite[Theorem 2.3]{at}   we deduce:
\begin{thm}[Thom isomorphism theorem] \label{ThomIsoReal} Suppose $X$ is a spin manifold of dimension $(8r-m)$. There is a natural isomorphism
\[
\varphi:\mathrm{KO}(X)=\mathrm{KR}(X)    \mapright{\sim} \mathrm{KR}^{m}(B(X),S(X)).
\]
\end{thm}
Now it is time to pass to index theory. 
It can be shown that the symbol of a real elliptic operator belongs to $\text{\rm KR}(BX,SX)$ where one complexifies the operator first; the involution covering $\iota$ comes then from complex conjugation. So, by construction, there is a forgetful map $\text{\rm KR}(BX,SX)\to\text{\rm K}(X)$ and the index theorem for complex bundles can be applied, but this gives nothing extra.
However, for families $X\to T$  the situation becomes different. The (analytic) index can be extended to a homomorphism
\[
\text{\rm ind} :  \text{\rm KR}(B(X/T),S(X/T))\to \text{\rm KO}(T) 
\]
covering the complex index map. But since the covering maps are in general not injective one gets extra information from the Index theorem for families of real elliptic operators \cite{AS4}.  It  states that  (complexified) analytic index  equals an  explicit expression in terms of Chern classes and  which can be called the topological index, $\text{\rm ind}_\tau $.

In the special case of a product family $X\times S^m \to S^m$ with $X$ a spin manifold of dimension $8r+m$, these two maps together with the above Thom isomorphism theorem induce a commutative diagram\footnote {Observe the change of sign in front of $m$.} 
\[
\begin{diagram}
\text{\rm KR} ((BX,SX) \times S^m)& \to  &\text{\rm KO}(S^m)\\
\downarrow&&\downarrow\\
\text{\rm KR}^{-m} (BX,SX)& \to &\text{\rm KO}^{-m}(\text{\rm point})\\
\Vline{}{}{2ex} && \Vline{}{}{2ex}\\
\text{\rm KO}^{ }(X)&\to&   \text{\rm KO}^{-m}(\text{\rm point}) .
\end{diagram} 
\]
In particular, if $d=8r+2$ and $m=2$ Bott periodicity gives two maps
\[
\text{\rm ind},\text{\rm ind}_\tau: \text{\rm KO}(X)   \to \text{\rm KO}^{-2}(\text{\rm point})=\Z/2\Z.
\]
These are equal and called the \emph{mod-$2$-index} for a family over $S^2$.
This  can in particular be applied   to real bundles  $E$; the  Dirac operator $\dirac D_E$ on $X$ with values in $E$ has  index $0$ (see  Cor.~\ref{Skew}) and a priori one does not expect information. But from such $E$ one can canonically construct  a family of Diracs depending on a complex parameter which then extends to the Riemann sphere $S^2$;  hence, the  above considerations  with $m=2$ apply.  It turns out (see \cite{AS4} for details) that the analytic index is the mod-$2$ dimension of the bundle  $\ker(\dirac D_{E})$ and the topological index  comes from the Gysin map  associated to $X\to $ point. Recall at this point    that for any  map  $f:X\to Y$ between compact \emph{spin} manifolds, there are    Gysin maps 
\[
f_!: \text{\rm KO}^*(X) \to \text{\rm KO}^{*-c}(Y), \quad c=\dim X-\dim Y.
\]
The upshot is
\begin{thm}\label{RealIndexThm}
Let $X$ be a compact spin manifold of dimension $2 \bmod 8$. Let $a_X:X\to{\text{\set{\rm point}}}$ be the constant map. For  $x\in \mathrm{KO}(X)$, let $j(x)$ be the mod--$2$ index of the Dirac operator with values in $x$.
Then there we have an equality of maps
\[
j= (a_X)_!: \mathrm{KO}(X)\to \mathrm{KO}^{-2}(\mathrm{point})=\Z/2.
\]
\end{thm}


\begin{thebibliography}{xxxxxxxxxx}

\bibitem[At]{at} Atiyah, M. F.:  $K$-theory and reality,  Quart. J. Math. Oxford Ser. (2)  \textbf{17}  367--386  (1966)
\bibitem[At1]{at0} Atiyah, M. F.:  \textit{$K$-theory},  W.A. Benjamin, Inc., New-York, Amsterdam (1967)
\bibitem[At2]{at2} Atiyah, M. F.: Bott periodicity, Quart.J. Math. \textbf{19} 113--140 (1968)
\bibitem[AH1]{ah1} Atiyah, M. F. and F. Hirzebruch: Riemann-Roch theorems for differentiable manifolds, Bull. AMS. \textbf{65} 276--281 (1959)
\bibitem[AH2]{ah2} Atiyah, M. F. and F. Hirzebruch: Vector bundles and homogeneous spaces, in \textsl{Differential Geometry}, Proc. Symp. Pure Math \textbf{3} Amer. Math. Soc., Providence R-I.  7--31 (1961)
\bibitem[AH3]{ah3} Atiyah, M. F. and F. Hirzebruch:  Charakterische Klassen und Anwendungen, Enseign. Math. II Ser. \textbf{7} 188--213 (1961)
\bibitem[AS1]{AS1} Atiyah, M. F. and I.M. Singer: The index of elliptic operators on compact manifolds. Bull. AMS. \textbf{69} 422--433 (1963)
\bibitem[AS2]{AS2} Atiyah, M. F.and I.M. Singer:  The index of elliptic operators. I. Ann. of Math. \textbf{87} 484--530 (1968) 
\bibitem[AS3]{AS3} Atiyah, M. F. and I.M. Singer: The index of elliptic operators IV Ann. Math. \textbf{93} 119-138 (1971)
\bibitem[AS4]{AS4} Atiyah, M. F. and I.M. Singer: The index of elliptic operators V, Ann. Math. \textbf{93} 139--149 (1971)
\bibitem[B-G-V]{heat} Berline, N., E. Getzler and M. Vergne: \textit{Heat Kernels and Dirac Operators}, Grundl. math. Wiss. \textbf{298}, Springer-Verlag,  Berlin etc. (1992)

\bibitem[Bor]{borcea} Borcea, C.: \textit{Calabi--Yau threefolds and complex multiplication}, in: Mirror Symmetry I (S.-T. Yau editor), 
Studies in Adv. Math. 9, 431--444 (1998)

\bibitem[Bott]{bo} Bott , R.: The stable homotopy group of the classicial groups, Ann. Math. \textbf{70} 313--337 (1959)

\bibitem[CSP]{periodbook}  Carlson, J., S. M\"uller-Stach  and C. Peters:  \textit{Period mappings and Period Domains}, Cambr. stud. in adv. math. \textbf{85} Cambr. Univ. press (2003)

\bibitem[Gr]{periodarticle} Griffiths, P.:{Periods of integrals on
algebraic manifolds, I, II }{Amer. J. Math.}{90}(1968){568--626,  805--865, respectively}



\bibitem[H-N-S]{gauge} Henningson, M,  B. E. W. Nilsson, and P. Salomonson:  Holomorphic Factorization  Of Correlation Functions In $(4k + 2)$-Dimensional (2k)-Form Gauge Theory   \texttt{hep-th/9908107} 

\bibitem[Hir]{hir} Hirzebruch, F.: \textit{Topological Methods in Algebraic Geometry}, Third Edition, Grundl. math. Wiss. \textbf{131}, Springer-Verlag,  Berlin etc. (1966)

\bibitem[McD-S]{symplec} McDuff, D. and D. Salamon: \textsl{Introduction to Symplectic Topology}, Oxford Math. Monogr. Clarendon Press, Oxford (1995)

\bibitem[Mo-Wi]{world} Moore, G. and E. Witten: Self-duality, Ramond-Ramond fields and $K$-theory, J. High Energy Phys. \textbf{5}, Paper 32, 32 pp. (2000)

\bibitem[PS]{mhb} Peters, C., J. Steenbrink:  \textit{Mixed Hodge Theory},   Ergebnisse Math., Springer Verlag, \textbf{52} (2008)


\bibitem[Mumford]{mum} Mumford, D.: \textsl{Abelian varieties}, Oxford University Press (1970)


\bibitem[Ru]{lazzeri} Rubei, E.:  {Lazzeri's  Jacobian of oriented compact riemannian manifolds 
Ark. Mat. \textbf{38} 381--397 (2000)} 

\bibitem[Te]{teichm} Teichm\"uller, O: Bestimmung der extremalen quasikonformen Abbildungen bei geschlossenen orientierten Riemannschen Fl\"achen, Abh. Preu\ss. Akad. Wiss., math.-naturw. Kl.4,  \textbf{4} 1--42 (1943) 

\bibitem[Warn]{warn} Warner, F.: \textit{Foundations  of Differentiable Manifolds and Lie Groups}, Graduate Texts in Math. \textbf{94},  Springer-Verlag,  Berlin etc. (1983)

\bibitem[Weil1]{weiljac} Weil, A.: On Picard varieties, Am. J. Math. \textbf{74} 865--894 (1962)

\bibitem[Weil2]{weil} Weil, A.: \textsl{Vari\'et\'es k\"ahleriennes}, Hermann, Paris (1958)

\bibitem[Witten]{wit} Witten, E.: Duality relations among topological effects in string theory. \texttt{hep-th/9912086}
\end{thebibliography}
\end{document}